\documentclass{article}
\usepackage{amsmath}
\usepackage{amssymb}
\usepackage{amsthm}
\usepackage{cite}
\usepackage[arrow,curve,matrix]{xy}
\begin{document}
\def\e#1\e{\begin{equation}#1\end{equation}}
\def\ea#1\ea{\begin{align}#1\end{align}}
\def\eq#1{{\rm(\ref{#1})}}
\theoremstyle{plain}
\newtheorem{thm}{Theorem}[section]
\newtheorem{lem}[thm]{Lemma}
\newtheorem{prop}[thm]{Proposition}
\newtheorem{cor}[thm]{Corollary}
\theoremstyle{definition}
\newtheorem{dfn}[thm]{Definition}
\def\dim{\mathop{\rm dim}\nolimits}
\def\Ker{\mathop{\rm Ker}}
\def\Image{\mathop{\rm Image}}
\def\red{{\rm red}}
\def\an{{\rm an}}
\def\na{{\rm na}}
\def\stk{{\rm stk}}
\def\orb{{\rm orb}}
\def\sign{{\rm sign}}
\def\Hom{\mathop{\rm Hom}\nolimits}
\def\Iso{\mathop{\rm Iso}\nolimits}
\def\Ad{\mathop{\rm Ad}}
\def\Stab{\mathop{\rm Stab}\nolimits}
\def\Spec{\mathop{\rm Spec}\nolimits}
\def\CF{\mathop{\rm CF}\nolimits}
\def\LCF{\mathop{\rm LCF}\nolimits}
\def\id{\mathop{\rm id}\nolimits}
\def\supp{\mathop{\rm supp}\nolimits}
\def\ge{\geqslant}
\def\le{\leqslant}
\def\Z{{\mathbin{\mathbb Z}}}
\def\Q{{\mathbin{\mathbb Q}}}
\def\C{{\mathbin{\mathbb C}}}
\def\K{{\mathbin{\mathbb K}}}
\def\CP{{\mathbin{\mathbb{CP}}}}
\def\KP{{\mathbin{\mathbb{KP}}}}
\def\D{{\mathbin{\mathfrak D}}}
\def\E{{\mathbin{\mathfrak E}}}
\def\F{{\mathbin{\mathfrak F}}}
\def\G{{\mathbin{\mathfrak G}}}
\def\H{{\mathbin{\mathfrak H}}}
\def\O{{\mathbin{\mathcal O}}}
\def\al{\alpha}
\def\be{\beta}
\def\ga{\gamma}
\def\de{\delta}
\def\io{\iota}
\def\la{\lambda}
\def\ka{\kappa}
\def\th{\theta}
\def\si{\sigma}
\def\De{\Delta}
\def\Ga{\Gamma}
\def\ts{\textstyle}
\def\sst{\scriptscriptstyle}
\def\sm{\setminus}
\def\ot{\otimes}
\def\iy{\infty}
\def\ra{\rightarrow}
\def\longra{\longrightarrow}
\def\t{\times}
\def\ci{\circ}
\def\ti{\tilde}
\def\bs{\boldsymbol}
\def\md#1{\vert #1 \vert}
\def\bmd#1{\big\vert #1 \big\vert}
\title{Constructible functions on Artin stacks}
\author{Dominic Joyce \\ Lincoln College, Oxford}
\date{}
\maketitle

\baselineskip 11.9pt plus .1pt

\section{Introduction}
\label{ca1}

Let $\K$ be an algebraically closed field, $X$ a $\K$-scheme, and
$X(\K)$ the set of closed points in $X$. A {\it constructible set\/}
$C\subseteq X(\K)$ is a finite union of subsets $Y(\K)$ for finite
type $\K$-subschemes $Y$ in $X$. A {\it constructible function}
$f:X(\K)\ra\Q$ has $f(X(\K))$ finite and $f^{-1}(c)$ constructible
for all $0\ne c\in f(X(\K))$. Write $\CF(X)$ for the $\Q$-vector
space of constructible functions on~$X$.

Let $\phi:X\ra Y$ and $\psi:Y\ra Z$ be morphisms of
$\C$-varieties. MacPherson \cite[Prop.~1]{MacP} defined a
$\Q$-linear {\it pushforward\/} $\CF(\phi):\CF(X)\ra\CF(Y)$
with $\bigl(\CF(\phi)\de_W\bigr)(y)=\chi_\an\bigl(\phi^{-1}(y)
\cap W\bigr)$ for subvarieties $W$ in $X$ and $y\in Y(\C)$,
where $\chi_\an$ is the topological Euler characteristic in
compactly-supported cohomology with the analytic topology,
and $\de_W$ the characteristic function of $W(\C)$ in $X(\C)$.
It satisfies $\CF(\psi\ci\phi)=\CF(\psi)\ci\CF(\phi)$, so
that $\CF$ is a {\it functor} from the category of
$\C$-varieties to the category of $\Q$-vector spaces.
This was extended to other fields $\K$ of characteristic
zero by Kennedy~\cite{Kenn}.

This paper generalizes these results to $\K$-{\it schemes} and {\it
algebraic $\K$-stacks} in the sense of Artin, for $\K$ of
characteristic zero. We introduce a notion of {\it pseudomorphism}
$\Phi$ between locally constructible sets in $\K$-schemes or
$\K$-stacks, generalizing morphisms. Pushforwards $\CF(\Phi)$ exist,
and pseudomorphisms seem very natural for constructible functions
problems.

The motivation for this is my series of papers
\cite{Joyc2,Joyc3,Joyc4,Joyc5}. Let coh$(P)$ be the abelian category
of {\it coherent sheaves} on a projective $\K$-scheme $P$, and
$(\tau,T,\le)$ a {\it stability condition} on coh$(P)$. Then the
moduli space ${\mathfrak{Obj}\kern .05em}_{{\rm coh}(P)}$ of sheaves
in coh$(P)$ is an Artin $\K$-stack, and the set ${\rm Obj}_{\rm
ss}^\al(\tau)$ of $\tau$-semistable sheaves in class $\al$ is a
constructible subset in ${\mathfrak{Obj}\kern .05em}_{{\rm
coh}(P)}$. We shall define invariants of $P,(\tau, T,\le)$ as
generalized Euler characteristics of ${\rm Obj}_{\rm ss}^\al(\tau)$,
and study identities they satisfy, and transformation laws under
change of stability condition.

To carry out this programme requires a theory of constructible sets
and functions in algebraic $\K$-stacks, and compatible notions of
Euler characteristic and pushforward. As I could not find these
tools in the literature, I develop them here. It seemed better to
write a stand-alone paper that others could use, rather than include
the material in the series~\cite{Joyc2,Joyc3,Joyc4,Joyc5}.

Section \ref{ca2} gives some background on schemes, varieties and
stacks. In \S\ref{ca3} we recall MacPherson's constructible
functions theory for $\C$-varieties, extend it to $\K$-schemes using
$l$-{\it adic cohomology} in place of cohomology with the analytic
topology for $\C$-varieties when $\K$ has characteristic zero, and
define and study {\it pseudomorphisms} between locally constructible
sets in $\K$-schemes. We also explain why the theory {\it cannot\/}
be extended to $\K$-schemes for $\K$ of positive characteristic.

Sections \ref{ca4} and \ref{ca5} extend these ideas to stacks.
An important difference between stacks and schemes is that
in an algebraic $\K$-stack $\F$ points $x\in\F(\K)$ have {\it
stabilizer groups} $\Iso_\K(x)$, which are algebraic $\K$-groups,
trivial if $\F$ is a $\K$-scheme. It turns out that there are
{\it many different ways} of including stabilizer groups when
extending Euler characteristics $\chi$ and pushforwards $\CF$
to stacks.

We highlight three interesting cases, the {\it
na\"\i ve pushforward\/} $\CF^\na$ which ignores
stabilizer groups, the {\it stack pushforward\/}
$\CF^\stk$ which is most natural in many stack
problems, and the {\it orbifold pushforward\/}
$\CF^\orb$, related to Deligne--Mumford stacks and
their crepant resolutions. Each is associated with a
notion of {\it Euler characteristic} $\chi^\na,\chi^\stk,
\chi^\orb$ of constructible sets in $\K$-stacks.

As $\chi^\stk,\CF^\stk$ involve weighting by $1/\chi(\Iso_\K(x))$,
the obvious definitions fail when $\chi(\Iso_\K(x))=0$. However, for
{\it representable} 1-morphisms $\phi:\F\!\ra\!\G$ we give a more
subtle definition of $\CF^\stk(\phi):\CF(\F)\!\ra\!\CF(\G)$ in
\S\ref{ca51}, which is always well-defined, and suffices for the
applications in \cite{Joyc3,Joyc4,Joyc5}. We also define {\it
pullbacks} $\psi^*$ by {\it finite type} 1-morphisms $\psi:\F\ra\G$,
and show pullbacks $\psi^*$ and pushforwards $\CF^\stk(\phi)$
commute in {\it Cartesian squares}.

A companion paper \cite{Joyc1} studies `stack functions' on Artin
stacks, which are a universal generalization of constructible
functions containing more information, and discusses how `motivic'
invariants of $\K$-varieties such as Euler characteristics and
virtual Poincar\'e polynomials are best extended to Artin stacks.

All $\K$-schemes and $\K$-stacks in this paper are assumed {\it
locally of finite type}.
\medskip

\noindent{\it Acknowledgements.} I thank Tom Bridgeland, Frances
Kirwan, Andrew Kresch, J\"org Sch\"urmann, Bertrand Toen and Burt
Totaro for useful conversations. I held an EPSRC Advanced Research
Fellowship whilst writing this paper.

\section{Schemes, varieties and stacks}
\label{ca2}

Fix an algebraically closed field $\K$ throughout. There are
four main classes of `spaces' over $\K$ used in algebraic
geometry, in increasing order of generality:
\begin{equation*}
\text{$\K$-varieties}\subset
\text{$\K$-schemes}\subset
\text{algebraic $\K$-spaces}\subset
\text{algebraic $\K$-stacks}.
\end{equation*}
Section \ref{ca21} gives a few definitions and facts on
$\K$-schemes and $\K$-varieties, and \S\ref{ca22} introduces
algebraic $\K$-stacks. Some good references for \S\ref{ca21}
are Hartshorne \cite{Hart}, and for \S\ref{ca22} are G\'omez
\cite{Gome} and Laumon and Moret-Bailly~\cite{LaMo}.

\subsection{Schemes and varieties}
\label{ca21}

We assume a good knowledge of $\K$-schemes and their morphisms,
following Hartshorne \cite{Hart}. We make the conventions that:
\begin{itemize}
\setlength{\itemsep}{0pt}
\setlength{\parsep}{0pt}
\item All $\K$-schemes in this paper are {\it locally of finite type}.
\item All $\K$-{\it subschemes} are {\it locally closed}, but not
necessarily closed.
\item A $\K$-{\it variety} is a reduced, irreducible, separated
$\K$-scheme of finite type.
\end{itemize}

\begin{dfn} For a $\K$-scheme $X$, write $X(\K)$ for
the set $\Hom(\Spec\K,X)$ of morphisms of $\K$-schemes
$\Spec\K\ra X$. Then $X(\K)$ is naturally identified
with the subset of {\it closed points} of the underlying
topological space of $X$. Elements of $X(\K)$ are also
called {\it geometric points} or $\K$-{\it points} of~$X$.

There is a natural identification $(X\t Y)(\K)\cong
X(\K)\t Y(\K)$. If $\phi:X\ra Y$ is a morphism of
$\K$-schemes, composition $\Spec\K\ra X\,{\smash{\buildrel
\phi\over\longra}}\,Y$ gives a natural map $\phi_*:
X(\K)\ra Y(\K)$. If $X$ is a $\K$-{\it subscheme} of
$Y$ then~$X(\K)\subseteq Y(\K)$.
\label{ca2def1}
\end{dfn}

Much of the paper will involve cutting schemes or stacks
into pieces. To do this we shall use two different notions
of {\it disjoint union}.

\begin{dfn} Let $X$ be a $\K$ scheme, and $\{X_i:i\in I\}$ a
family of $\K$-subschemes of $X$. We say that $X$ is the
{\it set-theoretic disjoint union} of the $X_i$
for $i\in I$ if~$X(\K)=\coprod_{i\in I}X_i(\K)$.

If $\{X_i:i\in I\}$ is a family of $\K$-schemes, we
define the {\it abstract disjoint union} of the
$X_i$ to be the $\K$-scheme $(X,\O_X)$, where $X$
is the disjoint union of the topological spaces $X_i$,
and $\O_X\vert_{X_i}=\O_{X_i}$. Then $X$ exists and is
unique up to isomorphism, and the $X_i$ are open and
closed $\K$-subschemes of $X$. Clearly, an abstract
disjoint union is a set-theoretic disjoint union, but
not necessarily vice versa. When we just say `disjoint
union' we mean set-theoretic disjoint union.
\label{ca2def2}
\end{dfn}

Here is a useful result of Rosenlicht \cite{Rose}, on the
existence of quotients of varieties by algebraic groups.

\begin{thm} Let\/ $\K$ be an algebraically closed field,
$X$ a $\K$-variety, $G$ an algebraic $\K$-group, and\/
$\rho:G\t X\ra X$ an algebraic action of\/ $G$ on $X$.
Then there exists a dense, Zariski open subset\/
$X'$ of\/ $X$, a $\K$-variety $Y$, and a surjective
morphism $\pi:X'\ra Y$ inducing a bijection between
$G$-orbits in $X'$ and\/ $\K$-points in $Y$, such that
any $G$-invariant rational function on $X'$ defined
at\/ $x\in X'$ is the pull-back of a rational function
on $Y$ defined at\/~$\pi(x)$.
\label{ca2thm1}
\end{thm}

\subsection{Algebraic stacks}
\label{ca22}

{\it Algebraic stacks} (also known as Artin stacks) were
introduced by Artin, generalizing {\it Deligne--Mumford stacks}.
For a good introduction to algebraic stacks see G\'omez
\cite{Gome}, and for a thorough treatment see Laumon and
Moret-Bailly \cite{LaMo}. As for schemes, we make the
convention that all algebraic $\K$-stacks in this paper
are {\it locally of finite type}, and $\K$-substacks are
{\it locally closed}.

Algebraic $\K$-stacks form a 2-{\it category}. That is, we have {\it
objects} which are $\K$-stacks $\F,\G$, and also two kinds of
morphisms, 1-{\it morphisms} $\phi,\psi:\F\ra\G$ between
$\K$-stacks, and 2-{\it morphisms} $A:\phi\ra\psi$ between
1-morphisms. An analogy to keep in mind is a 2-category of
categories, where objects are categories, 1-morphisms are functors
between the categories, and 2-morphisms are isomorphisms (natural
transformations) between functors.

We define the set of $\K$-{\it points} of a stack.

\begin{dfn} Let $\F$ be a $\K$-stack. Write $\F(\K)$ for the set of
2-isomorphism classes $[x]$ of 1-morphisms $x:\Spec\K\ra\F$.
Elements of $\F(\K)$ are called $\K$-{\it points}, or {\it geometric
points}, of $\F$. If $\phi:\F\ra\G$ is a 1-morphism then composition
with $\phi$ induces a map of sets~$\phi_*:\F(\K)\ra\G(\K)$.

For a 1-morphism $x:\Spec\K\ra\F$, the {\it stabilizer group}
$\Iso_\K(x)$ is the group of 2-morphisms $x\ra x$. When $\F$ is an
algebraic $\K$-stack, $\Iso_\K(x)$ is an {\it algebraic $\K$-group}.
We say that $\F$ {\it has affine geometric stabilizers} if
$\Iso_\K(x)$ is an affine algebraic $\K$-group for all 1-morphisms
$x:\Spec\K\ra\F$.

As an algebraic $\K$-group up to isomorphism, $\Iso_\K(x)$
depends only on the isomorphism class $[x]\in\F(\K)$ of $x$
in $\Hom(\Spec\K,\F)$. If $\phi:\F\ra\G$ is a 1-morphism,
composition induces a morphism of algebraic $\K$-groups
$\phi_*:\Iso_\K([x])\ra\Iso_\K\bigr(\phi_*([x])\bigr)$,
for~$[x]\in\F(\K)$.
\label{ca2def3}
\end{dfn}

One important difference in working with 2-categories rather
than ordinary categories is that in diagram-chasing one only
requires 1-morphisms to be 2-{\it isomorphic} rather than
{\it equal}. The simplest kind of {\it commutative diagram} is:
\begin{equation*}
\xymatrix@R=6pt{
& \G \ar@{=>}[d]^{\,F} \ar[dr]^\psi \\
\F \ar[ur]^\phi \ar[rr]_\chi && \H,
}
\end{equation*}
by which we mean that $\F,\G,\H$ are $\K$-stacks,
$\phi,\psi,\chi$ are 1-morphisms, and $F:\psi\ci\phi\ra\chi$
is a 2-isomorphism. Usually we omit $F$, and mean
that~$\psi\ci\phi\cong\chi$.

\begin{dfn} Let $\phi:\F\ra\H$, $\psi:\G\ra\H$ be 1-morphisms
of $\K$-stacks. Then one can define the {\it fibre product stack\/}
$\F\t_{\phi,\H,\psi}\G$, or $\F\t_\H\G$ for short, with 1-morphisms
$\pi_\F,\pi_\G$ fitting into a commutative diagram:
\e
\begin{gathered}
\xymatrix@R=-4pt{
& \F \ar[dr]^\phi \ar@{=>}[dd] \\
\F\t_\H\G
\ar[dr]_{\pi_\G} \ar[ur]^{\pi_\F} && \H.\\
& \G \ar[ur]_\psi \\
}
\end{gathered}
\label{ca2eq1}
\e
A commutative diagram
\begin{equation*}
\xymatrix@R=-4pt{
& \F \ar[dr]^\phi \ar@{=>}[dd] \\
\E
\ar[dr]_\eta \ar[ur]^\th && \H\\
& \G \ar[ur]_\psi \\
}
\end{equation*}
is a {\it Cartesian square} if it is isomorphic to \eq{ca2eq1}, so
there is a 1-isomorphism $\E\cong\F\t_\H\G$. Cartesian squares may
also be characterized by a universal property.
\label{ca2def4}
\end{dfn}

Here is a definition from Kresch \cite[Def.~3.5.3]{Kres},
slightly modified.

\begin{dfn} Let $\F$ be a finite type algebraic $\K$-stack,
and $\F^\red$ the associated reduced stack. We say that $\F$
{\it can be stratified by global quotient stacks} if
$\F^\red$ is the disjoint union of finitely many locally
closed substacks $U_i$ with each $U_i$ 1-isomorphic to a stack
of the form $[X_i/G_i]$, where $X_i$ is a $\K$-variety and
$G_i$ a smooth, connected, affine algebraic $\K$-group acting
linearly on $X_i$. For a stack to be the {\it disjoint union}
of a family of locally closed substacks is defined in
\cite[p.~22]{LaMo}. It implies that~$\F(\K)=\F^\red(\K)=
\coprod_iU_i(\K)$.
\label{ca2def5}
\end{dfn}

Kresch \cite[Prop.~3.5.2(ii)]{Kres} takes the $X_i$ to be
{\it quasiprojective schemes}, rather than varieties, but
this is equivalent to our definition. Kresch
\cite[Prop.~3.5.9]{Kres} characterizes stacks
stratified by global quotients.

\begin{thm} Let\/ $\F$ be a finite type algebraic $\K$-stack.
Then $\F$ can be stratified by global quotient stacks if
and only if\/ $\F$ has affine geometric stabilizers.
\label{ca2thm2}
\end{thm}

\section{Constructible functions on $\K$-schemes}
\label{ca3}

We now introduce {\it constructible sets\/} and {\it functions} on
$\K$-schemes, and the {\it pushforward\/} of constructible functions
by morphisms. Section \ref{ca31} defines (locally) constructible
sets and functions on $\K$-schemes. We explain the Euler
characteristic and pushforwards over $\C$ in \S\ref{ca32}, and over
other fields $\K$ in \S\ref{ca33}. Section \ref{ca34} defines {\it
pseudomorphisms}, a notion of morphism for (locally) constructible
sets, and pushforwards along pseudomorphisms.

Some references are Mumford \cite[p.~51]{Mumf} and Hartshorne
\cite[p.~94]{Hart} for constructible sets, and MacPherson
\cite{MacP}, Viro \cite{Viro} and Kennedy \cite{Kenn} for
constructible functions and the pushforward. As far as the
author can tell the ideas of \S\ref{ca33}--\S\ref{ca34} are
new, although elementary and probably obvious to experts.

\subsection{Constructible sets and functions on $\K$-schemes}
\label{ca31}

We define {\it constructible} and {\it locally constructible sets}.

\begin{dfn} Let $\K$ be an algebraically closed field,
and $X$ a $\K$-scheme. A subset $C\subseteq X(\K)$ is
called {\it constructible} if $C=\bigcup_{i\in I}X_i(\K)$,
where $\{X_i:i\in I\}$ is a finite collection of finite
type $\K$-subschemes $X_i$ of $X$. We call $S\subseteq X(\K)$
{\it locally constructible} if $S\cap C$ is constructible for
all constructible~$C\subseteq X(\K)$.
\label{ca3def1}
\end{dfn}

This is easily seen to be equivalent to a stronger definition,
where we take the union $C=\bigcup_{i\in I}X_i(\K)$ to be
{\it disjoint}, and the $X_i$ to be {\it separated}.

\begin{prop} Let\/ $X$ be a $\K$-scheme, and\/ $C\subseteq X(\K)$
a constructible subset. Then we may write $C=\coprod_{i\in I}X_i(\K)$,
where $\{X_i:i\in I\}$ is a finite collection of separated, finite
type $\K$-subschemes $X_i$ of\/~$X$.
\label{ca3prop1}
\end{prop}

The following properties of constructible sets in $\K$-{\it
varieties} are well known, \cite[p.~94]{Hart}, \cite[p.~51]{Mumf}.
Our extension to $\K$-schemes is straightforward.

\begin{prop} Let\/ $X,Y$ be $\K$-schemes, $\phi:X\ra Y$ a
morphism, and\/ $A,B\subseteq X(\K)$ be constructible subsets.
Then $A\cup B$, $A\cap B$ and\/ $A\sm B$ are constructible in
$X(\K)$, and\/ $\phi_*(A)$ is constructible in~$Y(\K)$.
\label{ca3prop2}
\end{prop}

Note that showing $\phi_*(A)$ constructible, and the stack
analogue in Proposition \ref{ca4prop2}, are the only places
we use the convention that $\K$-schemes and $\K$-stacks are
{\it locally of finite type}. Next we define ({\it locally})
{\it constructible functions}.

\begin{dfn} Let $X$ be a $\K$-scheme and $S\subseteq X(\K)$ be
locally constructible. A {\it constructible function} on $S$ is a
function $f:S\ra\Q$ such that $f(S)$ is finite and $f^{-1}(c)$ is a
constructible set in $S\subseteq X(\K)$ for each $c\in f(S)
\sm\{0\}$. Note that we do {\it not\/} require $f^{-1}(0)$ to be
constructible. Write $\CF(S)$ for the $\Q$-vector space of
constructible functions on $S$, and for brevity write $\CF(X)$
for~$\CF(X(\K))$.

A {\it locally constructible function} on $S$ is a function
$f:S\ra\Q$ such that $f\vert_C$ is constructible for all
constructible $C\subseteq S$. Equivalently, $f$ is locally
constructible if $f^{-1}(c)$ is locally constructible for all
$c\in\Q$, and $f(C)$ is finite for all constructible $C\subseteq S$.
Write $\LCF(S)$ for the $\Q$-vector space of locally constructible
functions on $S$, and $\LCF(X)$ for~$\LCF(X(\K))$.

Using Proposition \ref{ca3prop2} we see that products of (locally)
constructible functions are (locally) constructible, so $\CF(S)$ and
$\LCF(S)$ are {\it commutative $\Q$-algebras}, with $\CF(S)$ an {\it
ideal\/} in $\LCF(S)$. Note that $1\in\CF(S)$ if and only if $S$ is
constructible, so if it is not then $\CF(S)$ is an {\it algebra
without identity}.
\label{ca3def2}
\end{dfn}

Here are some remarks on this material:
\begin{itemize}
\setlength{\itemsep}{0pt}
\setlength{\parsep}{0pt}
\item To define constructible functions $f:X(\K)\ra\Q$ on
$\K$-schemes $X$ which are {\it not of finite type}, or
$f:S\ra\Q$ for $S$ {\it not constructible}, we must allow
$f^{-1}(0)$ to be non-constructible. If we did not there would
be {\it no} constructible functions on $X$ or $S$, not even~0.

For $X$ not of finite type we can think of $X(\K)$ as being
`large', or `unbounded'. Constructible functions $f:X(\K)\ra\Q$
are nonzero only on small, bounded subsets of $X(\K)$, and
$f^{-1}(0)$ is the remaining, large, unbounded part of~$X(\K)$.
\item We can also consider constructible functions with values
in $\Z$, or any other abelian group, ring or field. But for
simplicity we restrict to~$\Q$.
\end{itemize}

\subsection{Euler characteristics and pushforward for $\C$-schemes}
\label{ca32}

We define the {\it analytic Euler characteristic}~$\chi_\an$.

\begin{dfn} Let $X$ be a separated $\C$-scheme of finite type.
Then $X(\C)$ is a Hausdorff topological space with the {\it
analytic topology}. Write $\chi_\an(X)$ for the {\it Euler
characteristic} of $X(\C)$, in {\it compactly-supported cohomology}.
\label{ca3def3}
\end{dfn}

The following properties of $\chi_\an$ are well known.

\begin{prop} Let\/ $X,Y$ be separated $\C$-schemes of finite type. Then
\begin{itemize}
\setlength{\itemsep}{0pt}
\setlength{\parsep}{0pt}
\item[{\rm(i)}] If\/ $Z$ is a closed subscheme of\/ $X$ then
$\chi_\an(X)=\chi_\an(Z)+\chi_\an(X\sm Z)$.
\item[{\rm(ii)}] Suppose $X$ is the set-theoretic disjoint
union of subschemes $U_1,\ldots,U_m$. Then~$\chi_\an(X)=
\sum_{i=1}^m\chi_\an(U_m)$.
\item[{\rm(iii)}] $\chi_\an(X\t Y)=\chi_\an(X)\chi_\an(Y)$.
\item[{\rm(iv)}] If\/ $\phi:X\ra Y$ is a morphism which is
a locally trivial fibration in the analytic topology with
fibre $F$, then~$\chi_\an(X)=\chi_\an(F)\chi_\an(Y)$.
\item[{\rm(v)}] $\chi_\an(\C^m)=1$ and\/
$\chi_\an(\CP^m)=m+1$ for all\/~$m\ge 0$.
\end{itemize}
\label{ca3prop3}
\end{prop}

Now we can define {\it pushforwards\/} on $\C$-schemes.

\begin{dfn} Let $X$ be a $\C$-scheme and $C\subseteq X(\C)$ a
constructible subset. Proposition \ref{ca3prop1} gives
$C=\coprod_{i\in I}X_i(\C)$ for $\{X_i:i\in I\}$ finitely
many separated, finite type subschemes of $X$. Define
$\chi_\an(C)=\sum_{i\in I}\chi_\an(X_i)$. If $\{Y_j:j\in J\}$
is another choice from Proposition \ref{ca3prop1} then $X_i$
is the set-theoretic union of $X_i\cap Y_j$ for $j\in J$, so
Proposition \ref{ca3prop3}(ii) gives $\chi_\an(X_i)=\sum_{j\in J}
\chi_\an(X_i\cap Y_j)$. Hence
\begin{equation*}
\sum_{i\in I}\chi_\an(X_i)\!=\!\sum_{i\in I}\Bigl[\,\sum_{j\in J}
\chi_\an(X_i\cap Y_j)\Bigr]\!=\!\sum_{j\in J}\Bigl[\,\sum_{
i\in I}\chi_\an(X_i\cap Y_j)\Bigr]\!=\!\sum_{j\in J}\chi_\an(Y_j),
\end{equation*}
and $\chi_\an(C)$ is well-defined.

For $f\in\CF(X)$, define the {\it weighted Euler characteristic}
$\chi_\an(X,f)\in\Q$ by
\e
\chi_\an(X,f)=\ts\sum_{c\in f(X(\C))\sm\{0\}}c\,
\chi_\an\bigl(f^{-1}(c)\bigr).
\label{ca3eq1}
\e
This is well-defined as $f(X(\C))$ is finite and $f^{-1}(c)\subseteq
X(\C)$ is constructible for each $c\in f(X(\C))\sm\{0\}$. Clearly,
$f\mapsto\chi_\an(X,f)$ is a linear map~$\CF(X)\ra\Q$.

Now let $\phi:X\ra Y$ be a morphism of $\C$-schemes, and
$f\in\CF(X)$. Define the {\it pushforward\/} $\CF(\phi)f:
Y(\C)\ra\Q$ of $f$ to $Y$ by
\e
\CF(\phi)f(y)=\chi_\an\bigl(X,f\cdot\de_{\phi_*^{-1}(y)}\bigr)
\quad\text{for $y\in Y(\C)$.}
\label{ca3eq2}
\e
Here $\phi_*:X(\C)\ra Y(\C)$ is the induced map, $\phi_*^{-1}(y)
\subseteq X(\C)$ is the inverse image of $\{y\}$ under $\phi_*$,
and $\de_{\smash{\phi^{\smash{-1}}_*(y)}}$ is its characteristic
function. It is a locally constructible function, so
$f\cdot\de_{\smash{\phi^{\smash{-1}}_*(y)}}\in\CF(X)$, and
\eq{ca3eq2} is well-defined.
\label{ca3def4}
\end{dfn}

MacPherson \cite[Prop.~1]{MacP} gives an important property
of the pushforward for algebraic $\C$-varieties. The extension
to $\C$-schemes is straightforward. One can prove it by dividing
$X,Y$ into pieces upon which $\phi$ is a locally trivial
fibration in the analytic topology, and using
Proposition~\ref{ca3prop3}(ii),(iv).

\begin{thm} Let\/ $X,Y,Z$ be $\C$-schemes, $\phi:X\ra Y$
and\/ $\psi:Y\ra Z$ be morphisms, and\/ $f\in\CF(X)$. Then
$\CF(\phi)f$ is constructible, so $\CF(\phi):\CF(X)\ra\CF(Y)$
is a $\Q$-linear map. Also $\CF(\psi\ci\phi)=\CF(\psi)\ci\CF(\phi)$
as linear maps $\CF(X)\ra\CF(Z)$. Hence $\CF$ is a functor from the
category of\/ $\C$-schemes to the category of\/ $\Q$-vector spaces.
\label{ca3thm1}
\end{thm}

Viro \cite{Viro} gives an interesting point of view
on constructible functions. One can regard the Euler
characteristic as a {\it measure}, defined on
constructible sets. Then $\chi_\an(X,f)$ in \eq{ca3eq1}
is the {\it integral\/} of $f$ with respect to this
measure, and the pushforward $\CF(\phi)f$ integrates
$f$ over the fibres of~$\phi$.

\subsection{Extension to other fields $\K$}
\label{ca33}

To extend \S\ref{ca32} to other fields $\K$, we need a good notion
of Euler characteristic $\chi(X)$ for a separated $\K$-scheme $X$ of
finite type.

\begin{dfn} Let $\K$ be an algebraically closed field
of characteristic $p$, which may be zero, and fix a prime number
$l\ne p$. Write $\Q_l$ for the field of $l$-{\it adic rationals}.
Let $X$ be a separated $\K$-scheme of finite type. Then one may
define the {\it compactly-supported\/ $l$-adic cohomology groups}
$H^i_{\rm cs}(X,\Q_l)$ of $X$, for $i\ge 0$. The original reference
for \'etale and $l$-adic cohomology is Grothendieck et al.\
\cite{Grot}, and a good book is Milne \cite{Miln}. Define the {\it
Euler characteristic} $\chi(X)$ of $X$ to be
\e
\chi(X)=\ts\sum_{i=0}^{2\dim X}(-1)^i\dim_{\Q_l}H^i_{\rm cs}(X,\Q_l).
\label{ca3eq3}
\e
\label{ca3def5}
\end{dfn}

Here are some properties of $\chi$, generalizing
Proposition~\ref{ca3prop3}.

\begin{thm} Let\/ $\K$ be an algebraically closed field and
$X,Y$ be separated $\K$-schemes of finite type. Then
\begin{itemize}
\setlength{\itemsep}{0pt}
\setlength{\parsep}{0pt}
\item[{\rm(i)}] If\/ $Z$ is a closed subscheme of\/ $X$ then
$\chi(X)=\chi(Z)+\chi(X\sm Z)$.
\item[{\rm(ii)}] Suppose $X$ is the set-theoretic
disjoint union of subschemes $U_1,\ldots,U_m$.
Then~$\chi(X)=\sum_{i=1}^m\chi(U_i)$.
\item[{\rm(iii)}] $\chi(X\t Y)=\chi(X)\chi(Y)$.
\item[{\rm(iv)}] $\chi(X)$ is independent of the choice of\/ $l$
in Definition~\ref{ca3def5}.
\item[{\rm(v)}] When $\K=\C$ we have~$\chi(X)=\chi_\an(X)$.
\item[{\rm(vi)}] $\chi(\K^m)=1$ and\/ $\chi(\KP^m)=m+1$ for
all\/~$m\ge 0$.
\end{itemize}
\label{ca3thm2}
\end{thm}

\begin{proof} Part (i) comes from the long exact
sequence~\cite[4.XVII.5.1.16]{Grot}:
\begin{equation*}
\cdots\ra H^i_{\rm cs}(X\sm Z,\Q_l)\ra H^i_{\rm cs}(X,\Q_l)\ra
H^i_{\rm cs}(Z,\Q_l)\ra H^{i+1}_{\rm cs}(X\sm Z,\Q_l)\ra\cdots,
\end{equation*}
and (ii) follows from (i) by standard arguments. Part (iii) is a
consequence of the {\it K\"unneth formula}
\cite[4.XVII.5.4.3]{Grot}. Part (iv) is proved for $X$ {\it proper}
in \cite[5.VII.4.10]{Grot}. The general case follows from (i) as we
may write $X\cong\bar X\sm Z$ for $\bar X$ a proper separated
$\K$-scheme of finite type, and $Z$ a closed subscheme. Part (v)
follows from the {\it comparison theorem} \cite[4.XVI.4.1]{Grot}.
For (vi), calculation shows $H^i_{\rm cs}(\K^m,\Q_l)$ is $\Q_l$ if
$i=2m$ and 0 otherwise, so $\chi(\K^m)=1$. Then $\chi(\KP^m)=m+1$ by
(ii) and~$\KP^m=\coprod_{n=0}^m\K^n$.
\end{proof}

Here are the generalizations of Definition \ref{ca3def4}
and Theorem \ref{ca3thm1} to~$\K$.

\begin{dfn} Let $X$ be a $\K$-scheme and $C\subseteq X(\K)$ a
constructible subset. Write $C=\coprod_{i\in I}X_i(\K)$ as in
Proposition \ref{ca3prop1}, and define $\chi(C)=\sum_{i\in I}
\chi(X_i)$. This is well-defined as in Definition \ref{ca3def4},
using Theorem \ref{ca3thm2}(ii). For $f\in\CF(X)$, define the
{\it weighted Euler characteristic} $\chi(X,f)\in\Q$ by
\e
\chi(X,f)=\ts\sum_{c\in f(X(\K))\sm\{0\}}c\,\chi\bigl(f^{-1}(c)\bigr).
\label{ca3eq4}
\e
Then $f\mapsto\chi(X,f)$ is a linear map $\CF(X)\ra\Q$. Let
$\phi:X\ra Y$ be a morphism of $\K$-schemes. Define the
{\it pushforward\/} $\CF(\phi)f:Y(\K)\ra\Q$ of $f$ to $Y$ by
\e
\CF(\phi)f(y)=\chi\bigl(X,f\cdot\de_{\phi_*^{-1}(y)}\bigr)
\quad\text{for $y\in Y(\K)$.}
\label{ca3eq5}
\e
When $\K=\C$ these definitions agree with Definition
\ref{ca3def4} by Theorem~\ref{ca3thm2}(v).
\label{ca3def6}
\end{dfn}

There are several ways to prove the next theorem. One is to
use results of Kennedy \cite{Kenn}. He defines pushforwards
implicitly using intersections of Lagrangian cycles, but one
can show using base change and comparison theorems for $l$-adic
cohomology that his definition of $\CF(\phi)$ agrees with ours.
Another is to use Katz and Laumon \cite[Th.~3.1.2]{KaLa}, which
in characteristic zero relates pushforwards of constructible
sheaves and functions, so functoriality of $\CF$ follows from
that for sheaf pushforwards.

\begin{thm} Let\/ $\K$ be an algebraically closed field of
characteristic zero, $X,Y,Z$ be $\K$-schemes, $\phi:X\ra Y$
and\/ $\psi:Y\ra Z$ be morphisms, and\/ $f\in\CF(X)$. Then
$\CF(\phi)f$ is constructible, so $\CF(\phi):\CF(X)\ra\CF(Y)$ is
a $\Q$-linear map. Also $\CF(\psi\ci\phi)=\CF(\psi)\ci\CF(\phi)$,
so $\CF$ is a functor.
\label{ca3thm3}
\end{thm}

The last part of the theorem is {\it false} for $\K$ of
characteristic $p>0$, and I am grateful to J\"org Sch\"urmann
for the following explanation of why. The proof using
constructible sheaves fails because if $\mathcal L$ is a
locally constant $\Q_l$-sheaf of rank $r$ on a non-proper
$\K$-scheme $X$ for $l\ne p$, we need the fact that
\begin{equation*}
\chi\bigl(H^*_{\rm cs}(X,{\mathcal L})\bigr)
=r\cdot\chi\bigl(H^*_{\rm cs}(X,\Q_l)\bigr).
\end{equation*}
This holds in characteristic zero, but not in characteristic
$p>0$ without extra conditions on $\mathcal L$, which are
studied in Illusie~\cite{Illu}.

Here is a counterexample to Theorem \ref{ca3thm3} in positive
characteristic. Let $\K$ have characteristic $p\ge 2$, and
$\phi:\K\ra\K$ be the Artin--Schreier morphism $\phi:x\mapsto x^p-x$.
It is a $p$-fold \'etale covering of $\K$ by itself, so $\CF(\phi)1=p$
in $\CF(\K)$. Thus taking $\psi:\K\ra\Spec\K$ to be the projection
we have $\CF(\psi\ci\phi)1=1$ but $\CF(\psi)\ci\CF(\phi)1=p$ in
$\CF(\Spec\K)=\Q$, so~$\CF(\psi\ci\phi)\ne\CF(\psi)\ci\CF(\phi)$.

When $Z\!=\!\Spec\K$ and $\psi:Y\!\ra\!\Spec\K$ is the projection we
have $\CF(Z)\!=\!\Q$ and $\CF(\psi)g\!=\!\chi(Y,g)$ for $g\!\in\!
\CF(Y)$. So $\CF(\psi\ci\phi)\!=\!\CF(\psi)\!\ci\!\CF(\phi)$ gives a
relation between pushforwards and weighted Euler characteristics:

\begin{cor} Let\/ $\K$ have characteristic zero, $X,Y$ be
$\K$-schemes, $\phi:X\ra Y$ a morphism, and\/ $f\in\CF(X)$.
Then~$\chi(X,f)=\chi\bigl(Y,\CF(\phi)f\bigr)$.
\label{ca3cor1}
\end{cor}

\subsection{Extension to pseudomorphisms}
\label{ca34}

We define {\it pseudomorphisms}, a notion of morphism between
locally constructible sets that generalizes morphisms of schemes.

\begin{dfn} Suppose $\K$ is an algebraically closed field, $X,Y$
are $\K$-schemes and $S\subseteq X(\K)$, $T\subseteq Y(\K)$ are
locally constructible. Let $\Phi:S\ra T$ be a map, and
define the {\it graph\/} $\Ga_\Phi=\bigl\{(s,\Phi(s)):s\in
S\bigr\}$ in $X(\K)\t Y(\K)=(X\t Y)(\K)$. We call $\Phi$ a
{\it pseudomorphism} if $\Ga_\Phi\cap(C\t Y(\K))$ is
constructible for all constructible $C\subseteq X(\K)$.
This implies $\Ga_\Phi$ is locally constructible.

A pseudomorphism $\Phi$ is a {\it pseudoisomorphism}
if $\Phi$ is bijective and $\Phi^{-1}:T\ra S$ is a
pseudomorphism. When $S=X(\K)$ and $T=Y(\K)$ we
shall also call $\Phi:X\ra Y$ a pseudomorphism
(pseudoisomorphism) from $X$ to~$Y$.
\label{ca3def7}
\end{dfn}

When $X,Y$ are $\K$-{\it varieties}, pseudomorphisms $\Phi:X\ra Y$
coincide with {\it definable functions} in the model theory of
algebraic geometry. See for instance Marker \cite[\S 7.4]{Mark}, in
particular \cite[Lem.~7.4.7]{Mark} which shows that $\Phi$ equals a
{\it quasimorphism} of varieties on a nonempty open affine subset
$X_0$ of $X$, and a {\it morphism} if $\K$ has characteristic zero.
Here are some basic properties of pseudomorphisms. They are easily
proved using Proposition \ref{ca3prop2} and the projection morphisms
$X\t Y\ra Y$, $X\t Y\t Z\ra X\t Z$.

\begin{prop} Let\/ $\K$ be an algebraically closed field.
\begin{itemize}
\setlength{\itemsep}{0pt}
\setlength{\parsep}{0pt}
\item[{\rm(a)}] Let\/ $\phi:X\ra Y$ be a morphism (isomorphism)
of\/ $\K$-schemes. Then $\phi_*:X(\K)\ra Y(\K)$ is a
pseudomorphism (pseudoisomorphism).
\item[{\rm(b)}] Let\/ $X,Y$ be $\K$-schemes,
$S\subseteq X(\K)$, $T\subseteq Y(\K)$ be locally
constructible, $\Phi:S\ra T$ be a pseudomorphism, and\/
$C\subseteq S$ be constructible. Then $\Phi(C)$ is
constructible in $Y(\K)$. Also, if\/ $t\in T$ then
$C\cap\Phi^{-1}(t)$ is constructible in $X(\K)$. Hence,
$\Phi^{-1}(t)$ is locally constructible in~$X(\K)$.
\item[{\rm(c)}] Let\/ $X,Y,Z$ be $\K$-schemes,
$S\subseteq X(\K)$, $T\subseteq Y(\K)$,
$U\subseteq Z(\K)$ be locally constructible, and\/
$\Phi:S\ra T$, $\Psi:T\ra U$ be pseudo(iso)morphisms.
Then $\Psi\ci\Phi:S\ra U$ is a pseudo(iso)morphism.
\end{itemize}
\label{ca3prop4}
\end{prop}

We define {\it pushforwards} $\CF(\Phi):\CF(S)\ra\CF(T)$
along pseudomorphisms.

\begin{dfn} Let $X,Y$ be $\K$-schemes, $S\subseteq X(\K)$,
$T\subseteq Y(\K)$ be locally constructible, $\Phi:S\ra T$
a pseudomorphism, and $f\!\in\!\CF(S)$. Define the {\it
pushforward\/} $\CF(\Phi)f:T\ra\Q$ by
\e
\CF(\Phi)f(t)=\chi\bigl(S,f\cdot\de_{\Phi^{-1}(t)}\bigr)
\quad\text{for $t\in T$.}
\label{ca3eq6}
\e
Here $\de_{\smash{\Phi^{-1}(t)}}$ is the characteristic function
of $\Phi^{-1}(t)\subseteq S$ on $S$. By Proposition \ref{ca3prop4}(b)
$\de_{\smash{\Phi^{-1}(t)}}\in\LCF(S)$, and $f\in\CF(S)$, so
$f\cdot\de_{\Phi^{-1}(t)}\in\CF(S)$. Thus \eq{ca3eq6} is
well-defined, by Definition \ref{ca3def6}. If $\phi:X\ra Y$
is a morphism of $\K$-schemes then $\phi_*:X(\K)\ra Y(\K)$ is
a pseudomorphism by Proposition \ref{ca3prop4}(a), and
$\CF(\phi)$ in Definition \ref{ca3def6} coincides with
$\CF(\phi_*)$ above.
\label{ca3def8}
\end{dfn}

Here is the generalization of Theorems \ref{ca3thm1} and
\ref{ca3thm3} to pseudomorphisms.

\begin{thm} Let\/ $\K$ have characteristic zero, $X,Y,Z$ be
$\K$-schemes, $S\subseteq X(\K)$, $T\subseteq Y(\K)$,
$U\subseteq Z(\K)$ be locally constructible, $\Phi:S\ra T$,
$\Psi:T\ra U$ be pseudomorphisms, and\/ $f\in\CF(S)$. Then
$\CF(\Phi)f$ is constructible, so $\CF(\Phi):\CF(S)\ra\CF(T)$ is
$\Q$-linear, and\/~$\CF(\Psi\ci\Phi)\!=\!\CF(\Psi)\!\ci\!\CF(\Phi)$.
\label{ca3thm4}
\end{thm}

\begin{proof} Define $F_{\sst XY}:X(\K)\!\t\!Y(\K)\!\ra\!\Q$,
$F_{\sst XZ}:X(\K)\!\t\!Z(\K)\!\ra\!\Q$, $F_{\sst YZ}:Y(\K)
\!\t\!Z(\K)\!\ra\!\Q$ and $F_{\sst XYZ}:X(\K)\!\t\!Y(\K)\!\t
\!Z(\K)\!\ra\!\Q$ by
\ea
F_{\sst XY}(x,y)&=\begin{cases} f(x), & \text{$x\in S$ and
$y=\Phi(x)$,}\\ 0, & \text{otherwise,} \end{cases}
\label{ca3eq7}\\
F_{\sst XZ}(x,z)&=\begin{cases} f(x), & \text{$x\in S$ and
$z=\Psi\ci\Phi(x)$,}\\ 0, & \text{otherwise,} \end{cases}
\nonumber\\
F_{\sst YZ}(y,z)&=\begin{cases} (\CF(\Phi)f)(y), & \text{$y\in T$
and $z=\Psi(y)$,}\\ 0, & \text{otherwise.} \end{cases}
\nonumber\\
F_{\sst XYZ}(x,y,z)&=\begin{cases} f(x), & \text{$x\in S$,
$y=\Phi(x)$ and $z=\Psi(y)$}\\ 0, & \text{otherwise.} \end{cases}
\nonumber
\ea
Write $\Pi^{\sst XY}_{\sst Y}:X\t Y\!\ra\!Y$ for the
projection morphism, and so on. It is easy to show
$F_{\sst XY}\!\in\!\CF(X\t Y)$, so $\CF(\Pi^{\sst XY}_{\sst Y})
F_{\sst XY}\!\in\!\CF(Y)$ by Theorem \ref{ca3thm3}. But comparing
\eq{ca3eq5}--\eq{ca3eq7} shows $\CF(\Phi)f\!=\!\bigl(\CF(
\Pi^{\sst XY}_{\sst Y})F_{\sst XY}\bigr)\vert_T$. Therefore
$\CF(\Phi)f\!\in\!\CF(T)$, proving the first part.
For the second part, $F_{\sst XZ},F_{\sst YZ},F_{\sst XYZ}$
are constructible on $X\t Z,Y\t Z$ and $X\t Y\t Z$ in the
same way, and it is easy to prove that
\begin{align*}
\CF(\Psi\ci\Phi)f&=\!\bigl(\CF(\Pi^{\sst XZ}_{\sst Z})
F_{\sst XZ}\bigr)\vert_U, &
\CF(\Psi)\ci\CF(\Phi)f&=\!\bigl(\CF(\Pi^{\sst YZ}_{\sst Z})
F_{\sst YZ}\bigr)\vert_U,\\
\CF(\Pi^{\sst XYZ}_{\sst XZ})F_{\sst XYZ}&=F_{\sst XZ}
\quad\text{and} &
\CF(\Pi^{\sst XYZ}_{\sst YZ})F_{\sst XYZ}&=F_{\sst YZ}.
\end{align*}
But as $\Pi^{\sst XZ}_{\sst Z}\ci\Pi^{\sst XYZ}_{\sst XZ}=\Pi^{\sst
XYZ}_{\sst Z}=\Pi^{\sst YZ}_{\sst Z}\ci\Pi^{\sst XYZ}_{\sst YZ}$
Theorem \ref{ca3thm3} gives $\CF(\Pi^{\sst XZ}_{\sst Z})\ci
\CF(\Pi^{\sst XYZ}_{\sst XZ})=\CF(\Pi^{\sst XYZ}_{\sst Z})=
\CF(\Pi^{\sst YZ}_{\sst Z})\ci\CF(\Pi^{\sst XYZ}_{\sst YZ})$,
and the result follows.
\end{proof}

If $\Phi:S\ra T$ is a {\it pseudoisomorphism\/} then
$\Phi^{-1}(t)$ is a single point in \eq{ca3eq6}, giving
$\CF(\Phi)f(t)=f\ci\Phi^{-1}(t)$. We deduce:

\begin{cor} Let\/ $X,Y$ be $\K$-schemes, $S\subseteq X(\K)$,
$T\subseteq Y(\K)$ be locally constructible, and\/ $\Phi:S\ra T$
a pseudoisomorphism. Then $\CF(\Phi):\CF(S)\ra\CF(T)$ is an
isomorphism, with\/ $\CF(\Phi)f=f\ci\Phi^{-1}$
and\/~$\CF(\Phi^{-1})g=g\ci\Phi$.
\label{ca3cor2}
\end{cor}

The moral is that pseudoisomorphic (locally) constructible
sets are essentially {\it the same} from the point of view
of constructible functions. So in problems involving
constructible functions, we can work with ({\it locally})
{\it constructible sets up to pseudoisomorphism}, and
pseudomorphisms between them.

\section{Constructible functions on stacks}
\label{ca4}

We now generalize \S\ref{ca3} to stacks. Sections \ref{ca41}
and \ref{ca42} develop the basic definitions and properties
of constructible sets and functions, and show that any finite
type algebraic $\K$-stack $\F$ with affine geometric stabilizers
is pseudoisomorphic to a finite type $\K$-scheme. This enables
us to reduce to the scheme case of~\S\ref{ca3}.

An important difference between stacks and schemes is that points
$x\in\F(\K)$ in a $\K$-stack $\F$ have {\it stabilizer groups}
$\Iso_\K(x)$, which are trivial if $\F$ is a $\K$-scheme. There are
{\it many different ways} of including stabilizer groups when
extending Euler characteristics $\chi$ and pushforwards $\CF$ to
stacks. Section \ref{ca43} studies the simplest of these, the {\it
na\"\i ve} versions $\chi^\na, \CF^\na$, which just ignore
stabilizer groups. Given an {\it allowable weight function} $w$ upon
affine algebraic $\K$-groups, in \S\ref{ca44} we modify
$\chi^\na,\CF^\na$ to get $\chi_w,\CF_w$ by weighting by
$w_\F:x\!\mapsto\!w(\Iso_\K(x))$ on $\F(\K)$. Two special cases are
the {\it stack\/} versions $\chi^\stk,\CF^\stk$ which are most
natural in many problems, and the {\it orbifold\/} versions
$\chi^\orb,\CF^\orb$, related to Deligne--Mumford stacks and their
crepant resolutions.

\subsection{Basic definitions}
\label{ca41}

We begin by giving analogues for stacks of the major definitions
of~\S\ref{ca3}.

\begin{dfn} Let $\K$ be an algebraically closed field, and
$\F$ an algebraic $\K$-stack. We call $C\subseteq\F(\K)$
{\it constructible} if $C=\bigcup_{i\in I}\F_i(\K)$, where
$\{\F_i:i\in I\}$ is a finite collection of finite type
algebraic $\K$-substacks $\F_i$ of $\F$. We call $S\subseteq\F(\K)$
{\it locally constructible} if $S\cap C$ is constructible for all
constructible~$C\subseteq\F(\K)$.
\label{ca4def1}
\end{dfn}

Here is a partial analogue of Proposition \ref{ca3prop2},
proved in the same way.

\begin{lem} Let\/ $\F$ be an algebraic $\K$-stack and\/
$A,B\subseteq\F(\K)$ constructible subsets. Then $A\cup B$,
$A\cap B$ and\/ $A\sm B$ are constructible in~$\F(\K)$.
\label{ca4lem1}
\end{lem}

\begin{dfn} Let $\K$ be an algebraically closed field, $\F$ an
algebraic $\K$-stack, and $S\subseteq\F(\K)$ be locally
constructible. Call a function $f:S\ra\Q$ {\it constructible} if
$f(S)$ is finite and $f^{-1}(c)$ is a constructible set for each
$c\in f(S)\sm\{0\}$. Call $f:S\ra\Q$ {\it locally constructible} if
$f\vert_C$ is constructible for all constructible $C\subseteq
S\subseteq\F(\K)$. Write $\CF(S),\LCF(S)$ for the sets of (locally)
constructible functions on $S$. Using Lemma \ref{ca4lem1} we see
that $\CF(S), \LCF(S)$ are $\Q$-vector spaces. For brevity write
$\CF(\F), \LCF(\F)$ rather than~$\CF(\F(\K)),\LCF(\F(\K))$.

As in Definition \ref{ca3def2}, using Lemma \ref{ca4lem1} we see
that multiplication of functions makes $\CF(S),\LCF(S)$ into {\it
commutative $\Q$-algebras}, with $\CF(S)$ an ideal in $\LCF(S)$, and
$\CF(S)$ is an algebra {\it without identity} if $S$ is not
constructible.

Now let $\F,\G$ be algebraic $\K$-stacks, and $S\subseteq\F(\K)$,
$T\subseteq\G(\K)$ be locally constructible. Then $\F\t\G$ is an
algebraic $\K$-stack with $(\F\t\G)(\K)=\F(\K)\t\G(\K)$. Let
$\Phi:S\ra T$ be a map, and define the {\it graph\/} $\Ga_\Phi=
\bigl\{(s,\Phi(s)):s\in S\bigr\}$. We call $\Phi$ a {\it
pseudomorphism} if $\Ga_\Phi\cap(C\t\G(\K))$ is constructible
in $(\F\t\G)(\K)$ for all constructible $C\subseteq\F(\K)$. A
pseudomorphism $\Phi$ is a {\it pseudoisomorphism} if $\Phi$
is bijective and $\Phi^{-1}:T\ra S$ is a pseudomorphism.
\label{ca4def2}
\end{dfn}

These definitions agree with those of \S\ref{ca3} when
$\F,\G$ are $\K$-schemes.

\subsection{Constructible sets and pseudomorphisms in stacks}
\label{ca42}

We now extend properties of constructible sets and
pseudomorphisms in $\K$-schemes to algebraic $\K$-stacks
with affine geometric stabilizers.

\begin{prop} Let\/ $\F$ be a finite type algebraic $\K$-stack
with affine geometric stabilizers. Then there exist substacks
$\F_1,\ldots,\F_n$ of\/ $\F$ with\/ $\F(\K)=\coprod_{a=1}^n
\F_a(\K)$, $\K$-varieties $Y_1,\ldots,Y_n$, and\/ $1$-morphisms
$\phi_a:\F_a\ra Y_a$ with\/ $(\phi_a)_*:\F_a(\K)\ra Y_a(\K)$
bijective for~$a=1,\ldots,n$.
\label{ca4prop1}
\end{prop}

\begin{proof} By Theorem \ref{ca2thm2}, $\F$ can
be stratified by global quotient stacks. Thus
there exist finitely many substacks $U_i$ of $\F$
with $\F(\K)=\coprod_iU_i(\K)$, and each $U_i$
1-isomorphic to a quotient stack $[X_i/G_i]$, for
$X_i$ a $\K$-variety and $G_i$ an affine algebraic
$\K$-group acting on $X_i$. Theorem \ref{ca2thm1}
gives a dense open $G_i$-invariant $X_i'\subseteq X_i$,
and a morphism of $\K$-varieties $\pi_i:X_i'\ra Y_i$
inducing a bijection $X_i'(\K)/G_i\ra Y_i(\K)$.
We may write $X_i\sm X_i'$ as a disjoint union of
finitely many $G_i$-invariant $\K$-subvarieties
$X_{ij}$ with $\dim X_{ij}<\dim X_i$. Applying
Theorem \ref{ca2thm1} again gives a dense open
$G_i$-invariant $X_{ij}'\subseteq X_{ij}$ and
a morphism of $\K$-varieties $\pi_{ij}:X_{ij}'
\ra Y_{ij}$ inducing a bijection~$X_{ij}'(\K)/
G_{ij}\ra Y_{ij}(\K)$.

As the dimension decreases at each stage, this process
eventually yields finitely many substacks $\F_1,\ldots,\F_n$
of\/ $\F$ with\/ $\F(\K)=\coprod_{a=1}^n\F_a(\K)$, such that
$\F_a$ is 1-isomorphic to $[X_a/G_a]$ for $\K$-varieties
$X_a,Y_a$ and $G_a$ an algebraic $\K$-group acting on $X_a$,
and morphisms $\pi_a:X_a\ra Y_a$ inducing a bijection between
$X_a(\K)/G_a$ and $Y_a(\K)$. The 1-isomorphisms $\F_a\cong
[X_a/G_a]$ and $\pi_a$ combine to give a 1-morphism
$\phi_a:\F_a\ra Y_a$ with the properties we want.
\end{proof}

We extend the last part of Proposition \ref{ca3prop2} to stacks.

\begin{prop} Let\/ $\K$ be an algebraically closed field,
$\F,\G$ be algebraic $\K$-stacks with affine geometric stabilizers,
$\phi:\F\ra\G$ be a $1$-morphism, and\/ $C\subseteq\F(\K)$ be
constructible. Then $\phi_*(C)$ is constructible in~$\G(\K)$.
\label{ca4prop2}
\end{prop}

\begin{proof} By Definition \ref{ca4def1} $C=\bigcup_{i\in I}
\F_i(\K)$, where $\{\F_i:i\in I\}$ are finitely many finite
type substacks $\F_i$ of $\F$. So by Lemma \ref{ca4lem1} it
is enough to show each $\phi_*(\F_i(\K))$ is constructible.
As by convention $\G$ is locally of finite type it admits
an open cover $\{\G_j:j\in J\}$ of finite type substacks
$\G_j$. By Proposition \ref{ca4prop1}, for $a=1,\ldots,n_j$
there exist substacks $\G_{ja}$ of $\G_j$, $\K$-varieties
$Y_{ja}$ and 1-morphisms $\psi_{ja}:\G_{ja}\ra Y_{ja}$ with
$(\psi_{ja})_*$ bijective, such that~$\G_j(\K)=
\coprod_{a=1}^{n_j}\G_{ja}(\K)$.

Now $\bigl\{\phi^{-1}(\G_j):j\in J\bigr\}$ covers $\F_i$,
which is {\it quasicompact\/} as it is of finite type.
So there exists a finite subset $J_i\subseteq J$ such
that $\bigl\{\phi^{-1}(\G_j):j\in J_i\bigr\}$ covers $\F_i$.
Set $\F_{ija}=\F_i\cap\phi^{-1}(\G_{ja})$ for $j\in J_i$ and
$a=1,\ldots,n_j$. Then $\F_{ija}$ is a finite type $\K$-substack,
and $\phi_{ija}=\phi\vert_{\F_{ija}}:\F_{ija}\ra\G_{ja}$ a
1-morphism with $\phi_*(\F_i(\K))=\bigcup_{j\in J_i}\bigcup_{a=1}^{n_j}
(\phi_{ija})_*(\F_{ija}(\K))$. Hence by Lemma \ref{ca4lem1} it
suffices to show $(\phi_{ija})_*(\F_{ija}(\K))$ is constructible
in $\G_{ja}(\K)$ for all~$i,j,a$.

As $\F_{ija}$ is finite type it has an atlas $u_{ija}:U_{ija}
\ra\F_{ija}$ with $U_{ija}$ a finite type $\K$-scheme. Then
$\psi_{ja}\ci\phi_{ija}\ci u_{ija}:U_{ija}\ra Y_{ja}$ is a
morphism of $\K$-schemes. But $U_{ija}(\K)$ is constructible as
$U_{ija}$ is of finite type, so Proposition \ref{ca3prop2} shows
\begin{align*}
(\psi_{ja}\ci\phi_{ija}\ci u_{ija})_*\bigl(U_{ija}(\K)\bigr)&=
(\psi_{ja})_*\ci(\phi_{ija})_*\ci(u_{ija})_*\bigl(U_{ija}(\K)\bigr)\\
&=(\psi_{ja})_*\ci(\phi_{ija})_*\bigl(\F_{ija}(\K)\bigr)
\end{align*}
is constructible in $Y_{ja}(\K)$, where the second line follows
since $(u_{ija})_*$ is surjective as $u_{ija}$ is an atlas.
Now $\psi_{ja}:\G_{ja}\ra Y_{ja}$ is a finite type 1-morphism,
so it pulls back constructible subsets to constructible subsets.
Therefore
\begin{equation*}
(\psi_{ja})_*^{-1}\bigl((\psi_{ja})_*\ci(\phi_{ija})_*(\F_{ija}(\K))
\bigr)=(\phi_{ija})_*(\F_{ija}(\K))
\end{equation*}
is constructible in $\G_{ja}(\K)$, using $(\psi_{ja})_*$ a
bijection in the second step.
\end{proof}

Lemma \ref{ca4lem1} and Proposition \ref{ca4prop2} extend
Proposition \ref{ca3prop2} to algebraic $\K$-stacks with affine
geometric stabilizers. As the proof of Proposition \ref{ca3prop4}
depended only on Proposition \ref{ca3prop2}, it extends to
such stacks.

\begin{prop} Let\/ $\K$ be an algebraically closed field, and
$\F,\G,\H$ be algebraic $\K$-stacks with affine geometric stabilizers.
\begin{itemize}
\setlength{\itemsep}{0pt}
\setlength{\parsep}{0pt}
\item[{\rm(a)}] Let\/ $\phi:\F\ra\G$ be a $1$-morphism (\/$1$-isomorphism).
Then $\phi_*:\F(\K)\ra\G(\K)$ is a pseudomorphism (pseudoisomorphism).
\item[{\rm(b)}] Let\/ $S\subseteq\F(\K)$, $T\subseteq\G(\K)$ be
locally constructible, $\Phi:S\ra T$ be a pseudomorphism, and\/
$C\subseteq S$ be constructible. Then $\Phi(C)$ is constructible
in $\G(\K)$. Also, if\/ $t\in T$ then $C\cap\Phi^{-1}(t)$ is
constructible in $\F(\K)$. Hence, $\Phi^{-1}(t)$ is locally
constructible in~$\F(\K)$.
\item[{\rm(c)}] Let\/ $S\!\subseteq\!\F(\K)$, $T\!\subseteq\!\G(\K)$,
$U\!\subseteq\!\H(\K)$ be locally constructible, and\/ $\Phi:S\!\ra\!T$,
$\Psi:T\!\ra\!U$ be pseudo(iso)morphisms. Then $\Psi\!\ci\!\Phi:S\!\ra
\!U$ is a pseudo(iso)morphism.
\end{itemize}
\label{ca4prop3}
\end{prop}

The next proposition allows results about constructible sets
and functions on schemes to be easily extended to stacks.

\begin{prop} Let\/ $\K$ be an algebraically closed field,
$\F$ an algebraic $\K$-stack with affine geometric stabilizers,
and\/ $C\subseteq\F(\K)$ be constructible. Then $C$ is
pseudoisomorphic to $Y(\K)$ for a separated, finite type
$\K$-scheme~$Y$.
\label{ca4prop4}
\end{prop}

\begin{proof} Write $C=\coprod_{i\in I}\F_i(\K)$ for
$\F_i,i\in I$ finitely many finite type substacks in
$\F$. Proposition \ref{ca4prop1} gives substacks
$\F_{ia}$ in $\F_i$, $\K$-varieties $Y_{ia}$ and
1-morphisms $\phi_{ia}:\F_{ia}\ra Y_{ia}$ for
$a=1,\ldots,n_i$, with $\F_i(\K)=\coprod_{a=1}^{n_i}
\F_{ia}(\K)$, and $(\phi_{ia})_*$ bijective. Let $Y$ be
the abstract disjoint union of the $Y_{ia}$ for $i\in
I$ and $a=1,\ldots,n_i$, as in Definition \ref{ca2def2}.
It is a separated, finite type $\K$-scheme.
Define $\Phi:C\ra Y(\K)$ by $\Phi\vert_{\F_{ia}(\K)}=
(\phi_{ia})_*$ for all $i,a$. Then $\Phi$ is bijective,
as $(\phi_{ia})_*$ is. Proposition \ref{ca4prop3}(a)
shows $(\phi_{ia})_*$ is a pseudomorphism, so $\Phi$
is a pseudomorphism. As $\Phi$ is bijective and $C,Y(\K)$
constructible, $\Phi$ is a pseudoisomorphism.
\end{proof}

\subsection{The na\"\i ve Euler characteristic and pushforward}
\label{ca43}

Fix an algebraically closed field $\K$ of {\it characteristic zero}
for the rest of the section. We consider the simplest generalization
of $\chi(X),\chi(S,f),\CF(\Phi)$ to $\K$-stacks $\F$, which we call
{\it na\"\i ve} as it ignores the {\it stabilizer groups}
$\Iso_\K(x)$ for $x\in\F(\K)$. Here is the analogue of Definitions
\ref{ca3def6} and~\ref{ca3def8}.

\begin{dfn} Let $\F$ be an algebraic $\K$-stack with affine
geometric stabilizers, and $C\subseteq\F(\K)$ be constructible.
Then $C$ is pseudoisomorphic to $Y(\K)$ for a separated,
finite type $\K$-scheme $Y$ by Proposition \ref{ca4prop4}.
Define the {\it na\"\i ve Euler characteristic} $\chi^\na(C)$
by $\chi^\na(C)=\chi(Y)$, where $\chi(Y)$ is as in Definition
\ref{ca3def5}. If $Y'$ is another choice for $Y$ then $Y(\K)$
is pseudoisomorphic to $Y'(\K)$ by Proposition \ref{ca4prop3}(c),
so $\chi(Y)=\chi(Y')$. Thus $\chi^\na(C)$ is well-defined.
Now let $S\subseteq\F(\K)$ be locally constructible. For
$f\in\CF(S)$, define the {\it na\"\i ve weighted Euler
characteristic} $\chi^\na(S,f)\in\Q$ by
\begin{equation*}
\chi^\na(S,f)=\ts\sum_{c\in f(S)\sm\{0\}}c\,\chi^\na\bigl(f^{-1}(c)\bigr).
\end{equation*}

Suppose $\F,\G$ are algebraic $\K$-stacks with affine geometric stabilizers,
$S\subseteq\F(\K)$, $T\subseteq\G(\K)$ are locally constructible,
$\Phi:S\ra T$ is a pseudomorphism, and $f\in\CF(S)$. Define the
{\it na\"\i ve pushforward\/} $\CF^\na(\Phi)f:T\ra\Q$ of $f$ to
$T$ by
\e
\CF^\na(\Phi)f(t)=\chi^\na\bigl(S,f\cdot\de_{\Phi^{-1}(t)}\bigr)
\quad\text{for $t\in T$.}
\label{ca4eq1}
\e
Here $\de_{\smash{\Phi^{-1}(t)}}$ is the characteristic function
of $\Phi^{-1}(t)\subseteq S$ on $S$. As $\Phi^{-1}(t)$ is locally
constructible by Proposition \ref{ca4prop3}(b) we have
$\de_{\smash{\Phi^{-1}(t)}}\in\LCF(S)$, and $f\in\CF(S)$, so
$f\cdot\de_{\Phi^{-1}(t)}\in\CF(S)$. Thus \eq{ca4eq1} is
well-defined.
\label{ca4def3}
\end{dfn}

Here are the na\"\i ve generalizations of Theorem \ref{ca3thm4}
and Corollary~\ref{ca3cor1}.

\begin{thm} Let\/ $\F,\G,\H$ be algebraic $\K$-stacks
with affine geometric stabilizers, $S\!\subseteq\!\F(\K)$,
$T\!\subseteq\!\G(\K)$, $U\!\subseteq\!\H(\K)$ be locally
constructible, and\/ $\Phi:S\!\ra\!T$, $\Psi:T\!\ra\!U$ be
pseudomorphisms. If\/ $f\in\CF(S)$, then $\CF^\na(\Phi)f$
is a constructible function on $T$. Thus $\CF^\na(\Phi):
\CF(S)\ra\CF(T)$ is a $\Q$-linear map. Also $\CF^\na(\Psi
\ci\Phi)=\CF^\na(\Psi)\ci\CF^\na(\Phi)$ as linear
maps\/~$\CF(S)\ra\CF(U)$.
\label{ca4thm}
\end{thm}

\begin{proof} Define $A=\supp(f)\subseteq S$, $B=\Phi(A)\subseteq T$
and $C=\Psi(B)\subseteq U$. Then $A$ is constructible by Definition
\ref{ca4def2}, so $B,C$ are constructible by Proposition
\ref{ca4prop3}(b). By Proposition \ref{ca4prop4} there exist
separated, finite type $\K$-schemes $X,Y,Z$ and pseudoisomorphisms
$\al:A\ra X(\K)$, $\be:B\ra Y(\K)$ and $\ga:C\ra Z(\K)$. Then
$f\ci\al^{-1}\in\CF(X)$ as $f\vert_A\in\CF(A)$. By Proposition
\ref{ca4prop3}(c), $\be\ci\Phi\ci\al^{-1}:X(\K)\ra Y(\K)$ and
$\ga\ci\Psi\ci\be^{-1}:Y(\K)\ra Z(\K)$ are pseudomorphisms of
$\K$-schemes, so Theorem \ref{ca3thm4} gives $\CF(\be\ci\Phi\ci
\al^{-1})(f\ci\al^{-1})\in\CF(Y)$ and
\begin{equation*}
\CF(\ga\ci\Psi\ci\Phi\ci\al^{-1})(f\ci\al^{-1})=\CF(\ga
\ci\Psi\ci\be^{-1})\ci\CF(\be\ci\Phi\ci\al^{-1})(f\ci\al^{-1}).
\end{equation*}
Since $\be,\ga$ identify constructible sets and functions with
constructible sets and functions, these easily imply
$\bigl(\CF^\na(\Phi)f\bigr)\vert_B\in\CF(B)$ and
\begin{equation*}
\bigl(\CF^\na(\Psi\ci\Phi)f\bigr)\vert_C=\bigl(\CF^\na
(\Psi)\ci\CF^\na(\Phi)f\bigr)\vert_C.
\end{equation*}
As the unrestricted functions are zero outside $B,C$, the
theorem follows.
\end{proof}

\begin{cor} Let\/ $\F,\G$ be algebraic $\K$-stacks with affine
geometric stabilizers, $S\!\subseteq\!\F(\K)$, $T\!\subseteq\!
\G(\K)$ be locally constructible, $\Phi:S\!\ra\!T$ be a
pseudomorphism, and\/ $f\in\CF(S)$. Then~$\chi^\na(S,f)=
\chi^\na\bigl(T,\CF^\na(\Phi)f\bigr)$.
\label{ca4cor1}
\end{cor}

\subsection{Stabilizers $\Iso_\K(x)$ and weight functions}
\label{ca44}

We now discuss how to modify the na\"\i ve Euler characteristic
$\chi^\na$ and pushforward $\CF^\na$ of \S\ref{ca43} to take account
of {\it stabilizer groups} $\Iso_\K(x)$ for $x\in\F(\K)$. We do this
by inserting a {\it weight\/} $w_\F$ depending on $\Iso_\K(x)$. We
continue to fix $\K$ algebraically closed of characteristic zero.

\begin{dfn} Let $w:\{$affine algebraic $\K$-groups $G$$\}\ra\Q\cup
\{\iy\}$, $G\mapsto w(G)$ be a map with $w(G)=w(G')$ if $G\cong G'$.
If $\F$ is an algebraic $\K$-stack with affine geometric stabilizers,
define $w_\F:\F(\K)\ra\Q\cup\{\iy\}$ by $w_\F(x)=w(\Iso_\K(x))$.
We call $w$ an {\it allowable weight function} if $w_\F$ is a
{\it locally constructible function} on $\F$ with values in
$\Q\cup\{\iy\}$ for all $\F$. We also call $w$ {\it multiplicative}
if $w(G\t H)=w(G)w(H)$ for all affine algebraic $\K$-groups~$G,H$.
\label{ca4def4}
\end{dfn}

Here are the weighted analogues $\chi_w(C)$, $\chi_w(S,f)$,
$\CF_w(\Phi)f$ of $\chi^\na(C)$, $\chi^\na(S,f)$, $\CF^\na(\Phi)f$.
We allow $w$ to take the values $0,\iy$ to accommodate the examples
below. This means $\chi_w(C),\chi_w(S,f),\CF_w(\Phi)f$ are not
always defined.

\begin{dfn} Let $w$ be an allowable weight function, $\F,\G$
algebraic $\K$-stacks with affine geometric stabilizers,
$C\subseteq\F(\K)$ constructible, $S\subseteq\F(\K)$,
$T\subseteq\G(\K)$ locally constructible, and $\Phi:S\ra T$ a
pseudomorphism.

If $w_\F\ne\iy$ on $C$, define the $w$-{\it Euler characteristic}
$\chi_w(C)=\chi^\na(C,w_\F\vert_C)$. If $w_\F(c)=\iy$ for some
$c\in C$ we say $\chi_w(C)$ is {\it undefined}. For $f\in\CF(S)$
with $w_\F\ne\iy$ on $\supp f$, define the {\it weighted\/
$w$-Euler characteristic} $\chi_w(S,f)$ by $\chi_w(S,f)=
\chi^\na(S,w_\F f)=\chi^\na(\supp f,w_\F f)$, taking
$w_\F f=0$ outside $\supp f$ even where $w_\F=\iy$. If
$w_\F(s)=\iy$ for some $s\in\supp f$ we say $\chi_w(S,f)$
is {\it undefined}. If $w_\F\ne\iy$ on $S$ and $w_\G\ne 0$
on $T$ then $w_\F\in\LCF(S)$ and $w_\G^{-1}\in\LCF(T)$ by
Definition \ref{ca4def4}. Define
\e
\CF_w(\Phi)f=w_\G^{-1}\cdot\CF^\na(\Phi)(w_\F f)
\quad\text{for $f\in\CF(S)$.}
\label{ca4eq2}
\e
This is well-defined in $\CF(T)$ as $w_\F f\in\CF(S)$, so
$\CF^\na(\Phi)(w_\F f)\in\CF(T)$. Therefore $\CF_w(\Phi):
\CF(S)\ra\CF(T)$ is a $\Q$-linear map. If $w_\F(s)=\iy$
for some $s\in S$ or $w_\G(t)=0$ for some $t\in T$, we
say $\CF_w(\Phi)$ is {\it undefined}.
\label{ca4def5}
\end{dfn}

Then $\chi_w$ satisfies the following analogues of
Theorem~\ref{ca3thm2}(ii),(iii):

\begin{lem} Let\/ $w$ be an allowable weight function, and\/ $\F,\G$
algebraic $\K$-stacks with affine geometric stabilizers. Then
\begin{itemize}
\setlength{\itemsep}{0pt}
\setlength{\parsep}{0pt}
\item[{\rm(i)}] Suppose $C,D_1,\ldots,D_m\!\subseteq\!\F(\K)$ are
constructible with\/ $C=\coprod_{i=1}^mD_i$. Then $\chi_w(C)=
\sum_{i=1}^m\chi_w(D_i)$ if either side is defined.
\item[{\rm(ii)}] If\/ $w$ is multiplicative and\/
$C\subseteq\F(\K)$, $D\subseteq\G(\K)$ are constructible then
$\chi_w(C\t D)=\chi_w(C)\chi_w(D)$ if both sides are defined.
\end{itemize}
\label{ca4lem2}
\end{lem}

For the analogue of Theorem \ref{ca4thm}, from \eq{ca4eq2} we have
\begin{align*}
\CF_w(\Psi\ci\Phi)f&=w_\H^{-1}\cdot\CF^\na(\Psi\ci\Phi)(w_\F f)
=w_\H^{-1}\cdot\CF^\na(\Psi)\ci\CF^\na(\Phi)(w_\F f)\\
&=w_\H^{-1}\cdot\CF^\na(\Psi)\bigl[w_\G\cdot\CF_w(\Phi)f\bigr]
=\CF_w(\Psi)\ci\CF_w(\Phi)f
\end{align*}
by Theorem \ref{ca4thm}, provided everything is defined. So we deduce:

\begin{cor} Let\/ $w$ be an allowable weight function, $\F,\G,\H$
algebraic $\K$-stacks with affine geometric stabilizers,
$S\!\subseteq\!\F(\K)$, $T\!\subseteq\!\G(\K)$, $U\!\subseteq\!\H(\K)$
locally constructible with\/ $w_\F\!\ne\!\iy$ on $S$, $w_\G\!\ne\!0,\iy$
on $T$, $w_\H\!\ne\!0$ on $U$, and\/ $\Phi:S\!\ra\!T$,
$\Psi:T\!\ra\!U$ be pseudomorphisms. Then~$\CF_w(\Psi\ci\Phi)\!=\!
\CF_w(\Psi)\!\ci\!\CF_w(\Phi)$.
\label{ca4cor2}
\end{cor}

As in Corollary \ref{ca4cor1} we have:

\begin{cor} Let\/ $w$ be an allowable weight function, $\F,\G$
algebraic $\K$-stacks with affine geometric stabilizers,
$S\!\subseteq\!\F(\K)$, $T\!\subseteq\!\G(\K)$
locally constructible with\/ $w_\F\!\ne\!\iy$ on $S$,
$w_\G\!\ne\!0,\iy$ on $T$, $\Phi:S\ra T$ a
pseudomorphism, and\/ $f\in\CF(S)$. Then~$\chi_w(S,f)=
\chi_w\bigl(T,\CF_w(\Phi)f\bigr)$.
\label{ca4cor3}
\end{cor}

Here are two examples of multiplicative allowable weight functions.

\begin{prop}{\rm(a)} Define $e:\{$affine algebraic $\K$-groups$\}\ra\Z$
by $e(G)=\chi(G)$. Then $e$ is a multiplicative allowable weight function.
\vskip 3pt

\noindent{\rm(b)} If\/ $G$ is an affine algebraic $\K$-group, define
the adjoint action of\/ $G$ on itself by $\Ad(g)h=ghg^{-1}$. Then
the quotient\/ $[G/\Ad(G)]$ is an algebraic $\K$-stack of finite
type. Define $o(G)=\chi^\na\bigl([G/\Ad(G)]\bigr)$. Then $o:\{$affine
algebraic $\K$-groups$\}\ra\Z$ is a multiplicative allowable
weight function.
\label{ca4prop5}
\end{prop}

\begin{proof} Clearly $e$ and $o$ are well-defined and
multiplicative. Let $\F$ be an algebraic $\K$-stack with
affine geometric stabilizers. We must show $e_\F,o_\F\in
\LCF(\F)$, which holds provided $e_\G=e_\F\vert_{\G(\K)}$
and $o_\G=o_\F\vert_{\G(\K)}$ lie in $\CF(\G)$ for all finite
type $\K$-substacks $\G$ in $\F$. Theorem \ref{ca2thm2} gives
$\G(\K)=\coprod_{i\in I}U_i(\K)$, where $\{U_i:i\in I\}$ are
finitely many substacks of $\G$ with $U_i$ 1-isomorphic to
$[X_i/G_i]$ for $X_i$ a $\K$-variety and $G_i$ an affine
algebraic $\K$-group, acting on $X_i$ by~$\rho_i:X_i\t G_i\ra X_i$.

Write $\pi_i:X_i\ra U_i$ for the projection 1-morphism.
Let $Y_i$ be the inverse image under $\id_{X_i}\t\rho_i:
X_i\t G_i\ra X_i\t X_i$ of the diagonal in $X_i\t X_i$. Then
$Y_i$ is a finite type closed subscheme of $X_i\t G_i$. Let
$\si_i:Y_i\ra X_i$ be the restriction of the projection
$X_i\t G_i\ra X_i$. Then for each $x\in X_i(\K)$,
$\si_i^{-1}(x)=\{x\}\t\Stab_x(G_i)$, where $\Stab_x(G_i)$ is
the stabilizer subgroup of $x$ in~$G_i$.

Theorem \ref{ca3thm4} gives $\CF(\si_i)1\in\CF(X_i)$, as
$1\in\CF(Y_i)$. But for~$x\in X_i(\K)$
\begin{equation*}
\bigl(\CF(\si_i)1\bigr)(x)=\chi\bigl[\Stab_x(G_i)\bigr]=
\chi\bigl[\Iso_\K\bigl((\pi_i)_*(x)\bigr)\bigr]=e_\G
\bigl((\pi_i)_*(x)\bigr),
\end{equation*}
as $\si_i^{-1}(x)=\{x\}\t\Stab_x(G_i)$, and the stabilizer
group $\Iso_\K\bigl((\pi_i)_*(x)\bigr)$ in $U_i\cong[X_i/G_i]$
is $\Stab_x(G_i)$. Therefore $\CF(\pi_i)1=e_\G\ci(\pi_i)_*$
as maps $X_i(\K)\ra\Q$. Since $(\pi_i)_*$ is surjective this
implies that $e_\G\vert_{U_i(\K)}\in\CF(U_i)$, as
$(\pi_i)_*$ takes constructible sets to constructible sets
by Proposition \ref{ca4prop2}. But $\G(\K)=\coprod_{i\in I}U_i(\K)$
and $I$ is finite, so $e_\G\in\CF(\G)$. This proves~(a).

For (b), we form an algebraic $\K$-stack $\H_i$ with 1-morphisms
$Y_i\smash{{\buildrel\al_i\over\longra}\,\H_i\,{\buildrel\be_i
\over\longra}}\,X_i$ with $\pi_i\!=\!\be_i\!\ci\!\al_i$, such
that if $x\!\in\!X_i(\K)$ with $\Stab_x(G_i)\!=\!H$, so that
$\pi_i^{-1}(x)\!=\!\{x\}\t H$, then $\be_i^{-1}(x)\!=\!\{x\}
\!\t\![H/\Ad(H)]$, and $\al_i:\pi_i^{-1}(x)\!\ra\!\be_i^{-1}(x)$
is the projection $H\!\ra\![H/\Ad(H)]$. Then $\al_i$ is an atlas,
so $\H_i$ is of finite type. Thus $1\!\in\!\CF(\H_i)$, so
$\CF^\na(\be_i)1\!\in\!\CF(X_i)$ by Theorem \ref{ca4thm}.
But for~$x\!\in\!X_i(\K)$
\begin{align*}
\bigl(\CF^\na(\be_i)1\bigr)(x)&=\chi^\na\bigl(\bigl[
\Stab_x(G_i)/\Ad\bigl(\Stab_x(G_i)\bigr)\bigr]\bigr)\\
&=\chi^\na\bigl(\bigl[\Iso_\K\bigl((\pi_i)_*(x)\bigr)/
\Ad\bigl(\Iso_\K((\pi_i)_*(x))\bigr)\bigr]\bigr)
=o_\G\bigl((\pi_i)_*(x)\bigr).
\end{align*}
The rest of the proof is as for~(a).
\end{proof}

Other weight functions constructed from $e,o$ in a multiplicative
way are also multiplicative and allowable, such as $e^k,o^k,
\md{e}^k,\md{o}^k,\sign(e),\sign(o)$ and $e^ko^l$ for $k,l\in\Z$
with $k\,l>0$. We give special names to two interesting cases.

\begin{dfn} Let $w$ be an allowable weight function, $\F,\G$
algebraic $\K$-stacks with affine geometric stabilizers,
$C\!\subseteq\!\F(\K)$ constructible, $S\!\subseteq\!\F(\K)$,
$T\!\subseteq\!\G(\K)$ locally constructible, $\Phi:S\!\ra\!T$
a pseudomorphism, and~$f\!\in\!\CF(S)$.
\begin{itemize}
\setlength{\itemsep}{0pt}
\setlength{\parsep}{0pt}
\item[(a)] Define the {\it stack (weighted) Euler characteristic}
by $\chi^\stk(C)=\chi_{\smash{1/e}}(C)$ and $\chi^\stk(S,f)=
\chi_{1/e}(S,f)$, and the {\it stack pushforward\/} $\CF^\stk(\Phi)=
\CF_{1/e}(\Phi)$, in the sense of Definition \ref{ca4def4}, where $e$
is the weight function of Proposition \ref{ca4prop5}(a), and $1/e$ is
allowable. Here $\chi^\stk(C)$ is only defined if $\chi(\Iso_\K(c))\!
\ne\!0$ for all $c\!\in\!C$, and $\chi^\stk(S,f)$ only if
$\chi(\Iso_\K(s))\!\ne\!0$ for all $s\!\in\!\supp f$, and
$\CF^\stk(\Phi)$ only if $\chi(\Iso_\K(s))\!\ne\!0$ for all~$s\!\in\!S$.
\item[(b)] Define the {\it orbifold (weighted) Euler characteristic}
by $\chi^\orb(C)=\chi_o(C)$ and $\chi^\orb(S,f)=\chi_o(S,f)$, and
the {\it orbifold pushforward\/} $\CF^\orb(\Phi)=\CF_o(\Phi)$, where
$o$ is the allowable weight function of Proposition \ref{ca4prop5}(b).
As $o$ takes values in $\Z$, $\chi^\orb(C)$ and $\chi^\orb(S,f)$ are
always defined, and $\CF^\orb(\Phi)$ is defined if $o(\Iso_\K(t))\ne 0$
for all~$t\in T$.
\end{itemize}
\label{ca4def6}
\end{dfn}

The {\it stack Euler characteristic} $\chi^\stk$ and its pushforward
$\CF^\stk$ turn out to be the natural notions for the problems in
\cite{Joyc3,Joyc4,Joyc5}. If $X$ is a $\K$-variety and $G$ an
algebraic $\K$-group acting on $X$ with $\chi(G)\ne 0$, then
$\chi^\stk([X/G])=\chi(X)/\chi(G)$. It also has a universal property
in {\it Cartesian squares}, in~\S\ref{ca52}.

Unfortunately, as $\chi(G)=0$ for any algebraic $\K$-group $G$ with
$\K^\t=\K\sm\{0\}$ as a subgroup, $\chi^\stk(C), \chi^\stk(S,f)$ and
$\CF^\stk(\Phi)$ above are {\it undefined\/} in many interesting
situations, including everything in \cite{Joyc3,Joyc4,Joyc5}. But in
\S\ref{ca51} we will extend the definition of $\CF^\stk(\Phi)$ to
$\CF^\stk(\phi)$ for $\phi:\F\ra\G$ a {\it representable}
1-morphism, and this will be sufficient for the applications
of~\cite{Joyc3,Joyc4,Joyc5}.

For {\it Deligne--Mumford\/} stacks all stabilizer groups
are finite, and for $G$ finite $\chi(G)=\md{G}>0$, so that
$\chi^\stk,\CF^\stk$ are always defined. It is well-established
that for enumerative problems on Deligne--Mumford stacks one
counts a point $x\in\F(\K)$ with weight $1/\md{\Iso_\K(x)}$,
and $\chi^\stk$ generalizes this approach.

The {\it orbifold Euler characteristic} is the author's
attempt to generalize to stacks something already well
understood for orbifolds. Let $G$ be a finite group
acting on a compact manifold $M$, so that $M/G$ is
an {\it orbifold}. Dixon et al. \cite[p.~684]{DHVW}
observe the correct Euler characteristic of $M/G$ in
String Theory is
\e
\chi(M,G)=\frac{1}{\md{G}}\sum_{g,h\in G:gh=hg}\chi(M^{g,h}),
\label{ca4eq3}
\e
where~$M^{g,h}=\{x\in M:g\cdot x=h\cdot x=x\}$.

Atiyah and Segal \cite{AtSe} later interpreted $\chi(M,G)$
as the Euler characteristic of {\it equivariant K-theory}
$K_G(M)$. For a survey and further references on orbifold
Euler characteristics, see Roan \cite{Roan}. In particular,
it is believed and in many cases known that for a complex
orbifold $M/G$, $\chi(M,G)$ coincides with the Euler
characteristic $\chi(X)$ of any {\it crepant resolution}
of~$M/G$.

Let $M$ be a $\K$-scheme acted on by a finite group $G$.
Then $M^{g,h}$ is a subscheme of $M$, and \eq{ca4eq3}
makes sense. An easy calculation shows $\chi(M,G)=
\chi^\na([M/G],f)$, where $f\in\CF([M/G])$ is given by
\begin{align*}
f(x)&=\frac{\bmd{\bigl\{(g,h)\!\in\!\Iso_\K(x)^2:gh\!=\!hg\bigr\}}
}{\bmd{\Iso_\K(x)}}=\!\!\! \sum_{g\in\Iso_\K(x)}
\frac{\bmd{\{h\!\in\!\Iso_\K(x):gh\!=\!hg\}}}{\bmd{\Iso_\K(x)}}\\
&=\sum_{g\in\Iso_\K(x)}\frac{1}{\bmd{\Ad(\Iso_\K(x))g}}
=\bmd{\Iso_\K(x)/\Ad(\Iso_\K(x))}=o_{[M/G]}(x)
\end{align*}
for $x\in[M/G](\K)$. Hence $\chi(M,G)=\chi^\orb([M/G])$ by
Definitions \ref{ca4def4} and~\ref{ca4def6}.

Thus, our orbifold Euler characteristic $\chi^\orb([M/G])$
of the Deligne--Mumford stack $[M/G]$ agrees with the
physicists' orbifold Euler characteristic $\chi(M,G)$ of
the complex orbifold $M/G$ when $\K=\C$, but our notion
$\chi^{\rm orb}$ is also defined over other fields $\K$
and for more general stacks $\F$. It would be interesting
to know whether $\chi(M,G)$ being the Euler characteristic
of any crepant resolution over $\C$ extends using $\chi^\orb$
to other fields, or to more general stacks.

\section{Representable and finite type 1-morphisms}
\label{ca5}

Next we study stack pushforwards $\CF^\stk(\phi)$
by 1-{\it morphisms} $\phi:\F\!\ra\!\G$. Then $\phi_*:
\F(\K)\!\ra\!\G(\K)$ is a pseudomorphism, so the
obvious definition is $\CF^\stk(\phi)=\CF^\stk(\phi_*)$.
However, $\CF^\stk(\phi_*)$ is {\it undefined\/} if
$x\in\F(\K)$ with $\chi(\Iso_\K(x))=0$. Since $\chi(G)=0$ for
many affine algebraic $\K$-groups $G$, this is a serious drawback.
Instead, by using the {\it extra data} of the homomorphisms
$\phi_*:\Iso_\K(x)\ra\Iso_\K(\phi_*(x))$, in \S\ref{ca51} we
define $\CF^\stk(\phi)$ in many cases when $\CF^\stk(\phi_*)$
is undefined, in particular for all {\it representable}~$\phi$.

Section \ref{ca52} defines the {\it pullback\/}
$\psi^*:\CF(\G)\ra\CF(\F)$ for a {\it finite type} 1-morphism
$\psi:\F\ra\G$, and proves pullbacks $\psi^*$ and pushforwards
$\CF^\stk(\phi)$ commute in {\it Cartesian squares}. This will be an
important tool in \cite{Joyc3, Joyc4,Joyc5}. In \S\ref{ca53}, for
finite type $\phi:\F\ra\G$ we extend
$\CF^\na(\phi_*),\CF^\stk(\phi)$ to {\it locally constructible}
functions, with the usual functorial property.

Fix an algebraically closed field $\K$ of characteristic zero
for all of this section.

\subsection{Pushforwards by representable 1-morphisms}
\label{ca51}

Here is our definition of the stack pushforward $\CF^\stk(\phi)$
for a 1-morphism~$\phi$.

\begin{dfn} Let $\F,\G$ be algebraic $\K$-stacks with affine
geometric stabilizers and $\phi:\F\ra\G$ a 1-morphism. Then
for any $x\in\F(\K)$ we have a morphism $\phi_*:\Iso_\K(x)\ra
\Iso_\K(\phi_*(x))$ of affine algebraic $\K$-groups. The kernel
$\Ker\phi_*$ is an affine algebraic $\K$-group in $\Iso_\K(x)$,
so $\chi(\Ker\phi_*)$ is defined. The image $\phi_*(\Iso_\K(x))$
is an affine algebraic $\K$-group closed in $\Iso_\K(\phi_*(x))$,
so the quotient $\Iso_\K(\phi_*(x))/\phi_*(\Iso_\K(x))$ is a
quasiprojective $\K$-variety. Thus $\chi\bigl(\Iso_\K(
\phi_*(x))/\phi_*(\Iso_\K(x))\bigr)$ is also defined.

Suppose $\chi(\Ker\phi_*)\ne 0$ for all $x\in\F(\K)$.
Define $m_\phi:\F(\K)\ra\Q$ by
\e
m_\phi(x)=\frac{\chi\bigl(\Iso_\K(\phi_*(x))/
\phi_*(\Iso_\K(x))\bigr)}{\chi\bigl(\Ker\bigl(\phi_*:
\Iso_\K(x)\ra\Iso_\K(\phi_*(x))\bigr)\bigr)}
\quad\text{for $x\in \F(\K)$.}
\label{ca5eq1}
\e
An argument similar to Proposition \ref{ca4prop5} shows
$m_\phi\in\LCF(\F)$. Define the {\it stack pushforward\/}
$\CF^\stk(\phi):\CF(\F)\ra\CF(\G)$ by
\e
\CF^\stk(\phi)f=\CF^\na(\phi_*)(m_\phi\cdot f).
\label{ca5eq2}
\e
Here $m_\phi\cdot f\in\CF(\F)$ as $m_\phi\in\LCF(\F)$
and $f\in\CF(\F)$, so \eq{ca5eq2} is well-defined, and
$\CF^\stk(\phi):\CF(\F)\ra\CF(\G)$ is $\Q$-linear.
\label{ca5def1}
\end{dfn}

This agrees with the previous definition of $\CF^\stk(\phi_*)$
when it is defined, regarding $\phi_*:\F(\K)\ra\G(\K)$ as a
{\it pseudomorphism} by Proposition~\ref{ca4prop3}(a).

\begin{lem} In Definition \ref{ca5def1}, if\/ $\CF^\stk
(\phi_*)$ is defined in Definition \ref{ca4def6} then
$\CF^\stk(\phi):\CF(\F)\ra\CF(\G)$ is defined
and\/~$\CF^\stk(\phi_*)=\CF^\stk(\phi)$.
\label{ca5lem1}
\end{lem}

\begin{proof} Suppose $\CF^\stk(\phi_*)$ is defined. Then
$1/e_\F\!\ne\!\iy$ in $\F(\K)$, so $\chi(\Iso_\K(x))\ne 0$
for $x\!\in\!\F(\K)$. Now $\Ker\phi_*$ is normal in $\Iso_\K(x)$,
with quotient $\Iso_\K(x)/\Ker\phi_*$  naturally isomorphic
to $\phi_*(\Iso_\K(x))$. Hence
\begin{equation*}
\chi\bigl(\Iso_\K(x)\bigr)=\chi(\Ker\phi_*)\cdot
\chi\bigl(\phi_*(\Iso_\K(x))\bigr),
\end{equation*}
by general properties of $\chi$. As $\chi(\Iso_\K(x))\ne 0$
this implies $\chi(\Ker\phi_*)\ne 0$ for $x\in\F(\K)$, so
$\CF^\stk(\phi):\CF(\F)\ra\CF(\G)$ is defined. Similarly,
we have
\begin{equation*}
\chi\bigl(\Iso_\K(\phi_*(x))\bigr)=
\chi\bigl(\phi_*(\Iso_\K(x))\bigr)\cdot
\chi\bigl(\Iso_\K(\phi_*(x))/\phi_*(\Iso_\K(x))\bigr).
\end{equation*}

Dividing this equation by the previous one for $x\in\F(\K)$, which
is valid as $\chi(\Iso_\K(x))\ne 0$, and using \eq{ca5eq1} gives
\begin{equation*}
m_\phi(x)=\chi\bigl(\Iso_\K(\phi_*(x)\bigr)\big/
\chi\bigl(\Iso_\K(x)\bigr)=e_\G(\phi_*(x))/e_\F(x)
\quad\text{for $x\in\F(\K)$.}
\end{equation*}
It follows immediately from Definition \ref{ca4def6}
that~$\CF^\stk(\phi_*)=\CF^\stk(\phi)$.
\end{proof}

The functorial behaviour of Theorem \ref{ca4thm} holds
for~$\CF^\stk(\phi)$.

\begin{thm} Let\/ $\F,\G,\H$ be algebraic $\K$-stacks with affine
geometric stabilizers, and\/ $\phi:\F\!\ra\!\G$, $\psi:\G\!\ra\!\H$
$1$-morphisms. Suppose the kernels of\/ $\phi_*:\Iso_\K(x)
\!\ra\!\Iso_\K(\phi_*(x))$ for $x\!\in\!\F(\K)$ and\/
$\psi_*:\Iso_\K(y)\!\ra\!\Iso_\K(\psi_*(y))$ for $y\!\in\!\G(\K)$
have nonzero Euler characteristics. Then $\CF^\stk(\psi\!\ci\!\phi)
\!=\!\CF^\stk(\psi)\!\ci\!\CF^\stk(\phi)$ as well-defined linear
maps~$\CF(\F)\!\ra\!\CF(\H)$.
\label{ca5thm1}
\end{thm}

\begin{proof} Let $x\in\F(\K)$, and set $y=\phi_*(x)$ and
$z=\psi_*(y)$. Write
\begin{gather*}
G_x=\Iso_\K(x),\;\> G_y=\Iso_\K(y),\;\> G_z=\Iso_\K(z),\;\>
\phi_x=\phi_*:G_x\ra G_y,\\
\psi_y\!=\!\psi_*:G_y\!\ra\!G_z,\;\> K_{\phi,x}\!=\!\Ker(\phi_x),
\;\> K_{\psi,y}\!=\!\Ker(\psi_y),\;\>
K_{\psi\ci\phi,x}\!=\!\Ker(\psi_y\!\ci\!\phi_x),\\
I_{\phi,x}=\Image(\phi_x),\quad I_{\psi,y}=\Image(\psi_y),\quad
I_{\psi\ci\phi,x}=\Image(\psi_y\ci\phi_x).
\end{gather*}
Then $K_{\phi,x}$ is normal in $K_{\psi\ci\phi,x}$,
and the quotient $K_{\psi\ci\phi,x}/K_{\phi,x}$ is
isomorphic to $I_{\phi,x}\cap K_{\psi,y}$. So
general properties of $\chi$ give
\e
\chi(K_{\psi\ci\phi,x})=\chi(K_{\phi,x})\cdot
\chi(I_{\phi,x}\cap K_{\psi,y}).
\label{ca5eq3}
\e
The inclusions $I_{\phi,x}\cap K_{\psi,y}\subseteq
K_{\psi,y}$ and $I_{\psi\ci\phi,x}\subseteq I_{\psi,y}
\subseteq G_z$ imply that
\begin{align}
\chi(K_{\psi,y})&=\chi\bigl(K_{\psi,y}/(I_{\phi,x}\cap
K_{\psi,y})\bigr)\cdot\chi(I_{\phi,x}\cap K_{\psi,y})
\label{ca5eq4}\\
\text{and}\quad
\chi(G_z/I_{\psi\ci\phi,x})&=\chi(G_z/I_{\psi,y})\cdot
\chi(I_{\psi,y}/I_{\psi\ci\phi,x}).
\label{ca5eq5}
\end{align}

By assumption $\chi(K_{\phi,x}),\chi(K_{\psi,y})\ne 0$ for
$x\in\F(\K)$ and $y\in\G(\K)$, so $\CF^\stk(\phi):\CF(\F)
\ra\CF(\G)$ and $\CF^\stk(\psi):\CF(\G)\ra\CF(\H)$
are defined. As $\chi(K_{\psi,y})\ne 0$ equation \eq{ca5eq4}
gives $\chi(I_{\phi,x}\cap K_{\psi,y})\ne 0$, and this,
$\chi(K_{\phi,x})\ne 0$ and \eq{ca5eq3} show that
$\chi(K_{\psi\ci\phi,x})\ne 0$, which holds for all $x\in\F(\K)$.
Hence $\CF^\stk(\psi\ci\phi):\CF(\F)\ra\CF(\H)$ is defined.

Now $\psi_x^{-1}(I_{\psi\ci\phi,x})$ is an algebraic group with
$I_{\phi,x}\subseteq\psi_x^{-1}(I_{\psi\ci\phi,x})\subseteq G_y$, so
\begin{equation*}
\chi(G_y/I_{\phi,x})=\chi\bigl(G_y/\psi_x^{-1}(I_{\psi\ci\phi,x})\bigr)
\cdot\chi\bigl(\psi_x^{-1}(I_{\psi\ci\phi,x})/I_{\phi,x}\bigr).
\end{equation*}
But $\psi_x$ and $\ga(I_{\phi,x}\cap K_{\psi,y})\mapsto\ga\,I_{\phi,x}$
induce isomorphisms of homogeneous spaces
\begin{equation*}
G_y/\psi_x^{-1}(I_{\psi\ci\phi,x})\cong I_{\psi,y}/I_{\psi\ci\phi,x},
\quad K_{\psi,y}/(I_{\phi,x}\cap K_{\psi,y})\cong\psi_x^{-1}(I_{\psi
\ci\phi,x})/I_{\phi,x}.
\end{equation*}
Therefore the last two equations give
\e
\chi(G_y/I_{\phi,x})=\chi(I_{\psi,y}/I_{\psi\ci\phi,x})
\cdot\chi\bigl(K_{\psi,y}/(I_{\phi,x}\cap K_{\psi,y})\bigr).
\label{ca5eq6}
\e

Combining equations \eq{ca5eq1} and \eq{ca5eq3}--\eq{ca5eq6} yields
\begin{align*}
m_{\psi\ci\phi}(x)\!&=\!\frac{\chi(G_z/I_{\psi\ci\phi,x})}{\chi
(K_{\psi\ci\phi,x})}\!=\!\frac{\chi(G_z/I_{\psi,y})\chi(I_{\psi,y}/
I_{\psi\ci\phi,x})}{\chi(K_{\phi,x})\chi(I_{\phi,x}\cap K_{\psi,y})}
\cdot\frac{\chi\bigl(K_{\psi,y}/(I_{\phi,x}\!\cap\!K_{\psi,y})\bigr)
}{\chi\bigl(K_{\psi,y}/(I_{\phi,x}\!\cap\!K_{\psi,y})\bigr)}\\
&=\frac{\chi(G_z/I_{\psi,y})}{\chi(K_{\psi,y})}\cdot
\frac{\chi(G_y/I_{\phi,x})}{\chi(K_{\phi,x})}
=m_\psi(y)\cdot m_\phi(x).
\end{align*}
This identity is easily seen to be the extra ingredient
needed to modify the proof of Theorem \ref{ca4thm} to prove
that~$\CF^\stk(\psi\ci\phi)=\CF^\stk(\psi)\ci\CF^\stk(\phi)$.
\end{proof}

If $\phi:\F\ra\G$ is {\it representable} then $\phi_*:
\Iso_\K(x)\ra\Iso_\K(\phi_*(x))$ is an {\it injective}
morphism of algebraic $\K$-groups for all $x\in\F(\K)$.
Thus $\Ker\phi_*=\{1\}$, so $\chi(\Ker\phi_*)=1\ne 0$ for
all $x\in\F(\K)$, and $\CF^\stk(\phi)$ is defined. This gives:

\begin{thm} Let\/ $\F,\G,\H$ be algebraic $\K$-stacks with affine
geometric stabilizers, and\/ $\phi:\F\ra\G$, $\psi:\G\ra\H$
representable $1$-morphisms. Then $\psi\ci\phi:\F\ra\H$ is
representable, and\/ $\CF^\stk(\phi):\CF(\F)\ra\CF(\G)$,
$\CF^\stk(\psi):\CF(\G)\ra\CF(\H)$ and\/ $\CF^\stk(\psi\ci
\phi):\CF(\F)\ra\CF(\H)$ are well-defined linear maps
with\/~$\CF^\stk(\psi\ci\phi)=\CF^\stk(\psi)\ci\CF^\stk(\phi)$.
\label{ca5thm2}
\end{thm}

Also, for $\phi$ representable $m_\phi$ in \eq{ca5eq1}
takes values in $\Z$, so $\CF^\stk(\phi)$ maps $\Z$-valued
functions $\CF(\F)_\Z\!\subset\!\CF(\F)$ to $\Z$-valued
functions~$\CF(\G)_\Z\!\subset\!\CF(\G)$.

\subsection{Pullbacks by finite type 1-morphisms}
\label{ca52}

For {\it finite type} $\phi\!:\!\F\!\ra\!\G$ we can
pull back constructible functions from $\G$ to~$\F$.

\begin{dfn} Suppose $\phi:\F\ra\G$ is a finite type 1-morphism
of algebraic $\K$-stacks and $C\subseteq\G(\K)$ is constructible.
Then $C=\bigcup_{i\in I}\G_i(\K)$, where $\{\G_i:i\in I\}$ are
finitely many substacks of $\G$. Set $\F_i=\phi^*(\G_i)$, a
$\K$-substack of $\F$. Then $\F_i$ is of finite type, as $\phi,
\G_i$ are. Hence $\phi_*^{-1}(C)=\bigcup_{i\in I}\phi_*^{-1}
(\G_i(\K))=\bigcup_{i\in I}\F_i(\K)$ is constructible in $\F$. That
is, $\phi_*:\F(\K)\ra\G(\K)$ {\it pulls back constructible sets to
constructible sets}. Thus if $f\in\CF(\G)$ then $f\ci\phi_*$ lies in
$\CF(\F)$. Define the {\it pullback\/} $\phi^*:\CF(\G) \ra\CF(\F)$
by $\phi^*(f)=f\ci\phi_*$. Pullbacks commute with multiplication of
functions, that is, $\phi^*(fg)=\phi^*(f)\phi^*(g)$. If also
$\psi:\G\ra\H$ is a finite type 1-morphism, it is immediate
that~$(\psi\!\ci\!\phi)^*\!=\!\phi^*\!\ci\!\psi^*:
\CF(\H)\!\ra\!\CF(\F)$.
\label{ca5def2}
\end{dfn}

It is an interesting question how pullbacks $\psi^*$ and
pushforwards $\CF^\stk(\phi)$ are related. The next theorem shows
they commute in {\it Cartesian squares}, as in Definition
\ref{ca2def4}. It will be an important tool in
\cite{Joyc3,Joyc4,Joyc5}. The theorem would {\it not\/} hold if we
replaced $\CF^\stk(\eta),\CF^\stk(\phi)$ in \eq{ca5eq7} by
$\CF^\na(\eta),\CF^\na(\phi)$, or pushforwards defined using some
other weight function. This supports our claim that $\CF^\stk$ is
the most natural pushforward in many stack problems.

\begin{thm} Let\/ $\E,\F,\G,\H$ be algebraic $\K$-stacks with
affine geometric stabilizers. If
\e
\begin{gathered}
\xymatrix{
\E \ar[r]_\eta \ar[d]^\th & \G \ar[d]_\psi \\
\F \ar[r]^\phi & \H
}
\end{gathered}
\quad
\begin{gathered}
\text{is a Cartesian square with}\\
\text{$\eta,\phi$ representable and}\\
\text{$\th,\psi$ of finite type, then}\\
\text{the following commutes:}
\end{gathered}
\quad
\begin{gathered}
\xymatrix@C=35pt{
\CF(\E) \ar[r]_{\CF^\stk(\eta)} & \CF(\G) \\
\CF(\F) \ar[r]^{\CF^\stk(\phi)} \ar[u]_{\th^*}
& \CF(\H). \ar[u]^{\psi^*}
}
\end{gathered}
\label{ca5eq7}
\e
\label{ca5thm3}
\end{thm}

\begin{proof} Let $C\subseteq\F(\K)$ be constructible, and
$\de_C\in\CF(\F)$ be its characteristic function. We shall
prove that
\e
\bigl(\CF^\stk(\eta)\ci\th^*\bigr)\de_C=
\bigl(\psi^*\ci\CF^\stk(\phi)\bigr)\de_C.
\label{ca5eq8}
\e
As $\CF^\stk(\eta)\ci\th^*$, $\psi^*\ci\CF^\stk(\phi)$ are
linear and such $\de_C$ generate $\CF(\F)$, this implies
$\CF^\stk(\eta)\ci\th^*=\psi^*\ci\CF^\stk(\phi)$, as we want.

Define $B=\th_*^{-1}(C)$. Since $\th$ is of finite type
$B$ is constructible, as in Definition \ref{ca5def2},
and $\th^*(\de_C)=\de_B$, the characteristic function of $B$.
Let $x\in\G(\K)$, and define $y=\psi_*(x)$, $B_x=B\cap
\eta_*^{-1}(\{x\})$ and $C_y=C\cap\phi_*^{-1}(\{y\})$. Then
$B_x,C_y$ are constructible, as $B,C$ are and $\eta_*^{-1}
(\{x\}),\phi_*^{-1}(\{y\})$ are locally constructible.
Write $\de_{B_x},\de_{C_y},\de_{\smash{
\eta_*^{-1}(x)}},\de_{\smash{\phi_*^{-1}(y)}}$ for the
characteristic functions. Then
\e
\begin{gathered}
\bigl((\CF^\stk(\eta)\ci\th^*)\de_C\bigr)(x)
=\bigl(\CF^\stk(\eta)(\de_C\ci\th_*)\bigr)(x)
=\bigl(\CF^\stk(\eta)\de_B\bigr)(x)=\\
\bigl(\CF^\na(\eta_*)(m_\eta\!\cdot\!\de_B)\bigr)(x)
\!=\!\chi^\na\bigl(\E,m_\eta\!\cdot\!\de_{B_x}\bigr)\!=\!
\chi^\na\bigl(\F,\CF^\na(\th)(m_\eta\!\cdot\!\de_{B_x})\bigr),
\end{gathered}
\label{ca5eq9}
\e
by \eq{ca4eq1}, \eq{ca5eq2} and Corollary \ref{ca4cor1}, where
$m_\eta$ is defined in \eq{ca5eq1}. Similarly we have
\e
\bigl((\psi^*\ci\CF^\stk(\phi))\de_C\bigr)(x)=
\chi^\na\bigl(\F,m_\phi\cdot\de_{C_y}\bigr).
\label{ca5eq10}
\e

We shall prove that
\e
\CF^\na(\th)(m_\eta\cdot\de_{B_x})=m_\phi\cdot\de_{C_y}
\quad\text{in $\CF(\F)$.}
\label{ca5eq11}
\e
If $z\in\F(\K)\sm C_y$ then both sides of \eq{ca5eq11} are zero
at $z$. So let $z\in C_y$. Then $\th_*^{-1}(\{z\})\cap B_x\!=\!
\th_*^{-1}(\{z\})\!\cap\!\eta_*^{-1}(\{x\})$, so by \eq{ca4eq1}
equation \eq{ca5eq11} at $z$ reduces to
\e
\chi^\na\bigl(\E,m_\eta\cdot\de_{\eta_*^{-1}(x)}\cdot
\de_{\th_*^{-1}(z)}\bigr)=m_\phi(z).
\label{ca5eq12}
\e

Define $G_x=\Iso_\K(x)$, $G_y=\Iso_\K(y)$ and $G_z=\Iso_\K(z)$,
as algebraic $\K$-groups. Since $\psi_*(x)=y$ and $\phi_*(z)=y$
we have homomorphisms $\psi_*:\Iso_\K(x)\ra\Iso_\K(y)$ and
$\phi_*:\Iso_\K(z)\ra\Iso_\K(y)$. Write these as $\psi_x:G_x
\ra G_y$ and $\phi_z:G_z\ra G_y$. Then $\phi_z$ is injective,
as $\phi$ is representable, so $\chi(\Ker\phi_z)=\{1\}$ and
\eq{ca5eq1} gives
\e
m_\phi(z)=\chi\bigl(G_y/\phi_z(G_z)\bigr).
\label{ca5eq13}
\e
As \eq{ca5eq7} is Cartesian $\E$ is 1-isomorphic to $\F\t_\H\G$
by Definition \ref{ca2def4}. By definition of fibre products
we find $\eta_*^{-1}(\{x\})\cap\th_*^{-1}(\{z\})$ is naturally
isomorphic to $\psi_x(G_x)\backslash G_y/\phi_z(G_z)$, a {\it
biquotient}. The stabilizer groups are given by
\begin{equation*}
\Iso_\K\bigl(\psi_x(G_x)\be\phi_z(G_z)\bigr)=\bigl\{(\al,\ga)\in
G_x\t G_y:\psi_x(\al)\be=\be\phi_z(\ga)\bigr\}
\quad\text{for $\be\in G_y$,}
\end{equation*}
and the group homomorphism $\eta_*:\Iso_\K\bigl(\psi_x(G_x)\be
\phi_z(G_z)\bigr)\ra\Iso_\K(x)=G_x$ is given by $(\al,\ga)\mapsto\al$.
It is injective as $\phi_z$ is injective. Thus \eq{ca5eq1} yields
\e
m_\eta\bigl(\psi_x(G_x)\be\phi_z(G_z)\bigr)=
\chi\bigl(G_x/\{\al\in G_x:\psi_x(\al)\be\phi_z(G_z)=
\be\phi_z(G_z)\}\bigr).
\label{ca5eq14}
\e

Let $\Pi_{x,y,z}:G_y/\phi_z(G_z)\ra\psi_x(G_x)\backslash
G_y/\phi_z(G_z)$ be the natural projection. Then the fibre of
$\Pi_{x,y,z}$ over $\psi_x(G_x)\be\phi_z(G_z)$ is isomorphic
to $G_x/\{\al\in G_x:\psi_x(\al)\be\phi_z(G_z)=\be\phi_z(G_z)\}$.
So \eq{ca5eq14} implies that $\CF^\na(\Pi_{x,y,z})1=m_\eta$ in
$\CF\bigl(\psi_x(G_x)\backslash G_y/\phi_z(G_z)\bigr)$. Therefore
\begin{gather*}
\chi^\na\bigl(\E,m_\eta\cdot\de_{\eta_*^{-1}(x)}\cdot
\de_{\th_*^{-1}(z)}\bigr)=
\chi^\na\bigl(\psi_x(G_x)\backslash G_y/\phi_z(G_z),m_\eta\bigr)=\\
\chi^\na\bigl(\psi_x(G_x)\backslash G_y/\phi_z(G_z),\CF^\na(
\Pi_{x,y,z})1\bigr)=\chi^\na\bigl(G_y/\phi_z(G_z),1\bigr)=m_\phi(z),
\end{gather*}
by Corollary \ref{ca4cor1} and \eq{ca5eq13}. This proves \eq{ca5eq12},
and hence \eq{ca5eq11}. Equations \eq{ca5eq9}--\eq{ca5eq11} then
give \eq{ca5eq8} at $x$, as we have to prove.
\end{proof}

Note that by \cite[Rem.~4.14.1 \& Lem.~3.11 \&
Rem.~4.17(2)]{LaMo}, in a Cartesian square \eq{ca5eq7} of
algebraic $\K$-stacks, if $\phi$ is representable then $\eta$
is representable, and if $\psi$ is of finite type then $\th$
is of finite type. Thus it is enough to suppose only that
$\phi$ is representable and $\psi$ of finite type
in~\eq{ca5eq7}.

\subsection{Pushforwards of locally constructible functions}
\label{ca53}

Next we observe that if $\phi:\F\ra\G$ is of {\it finite type} then
the definitions of $\CF^\na(\phi_*)f,\CF^\stk(\phi)f$ in
\eq{ca4eq1}, \eq{ca5eq2} make sense for $f$ only {\it locally}
constructible.

\begin{dfn} Let $\phi:\F\!\ra\!\G$ be a finite type 1-morphism
of algebraic $\K$-stacks with affine geometric stabilizers. For
$f\!\in\!\LCF(\F)$, define $\LCF^\na(\phi)f$~by
\e
\LCF^\na(\phi)f(x)=\chi^\na\bigl(\F,f\cdot\de_{\phi_*^{-1}(x)}\bigr)
\quad\text{for $x\in\F(\K)$,}
\label{ca5eq15}
\e
following \eq{ca4eq1}. This is well-defined as $\phi_*^{-1}(\{x\})$
is constructible since $\phi$ is of finite type. Thus
$\de_{\smash{\phi_*^{-1}(x)}}\in\CF(\F)$ and $f\in\LCF(\F)$,
giving $f\cdot\de_{\smash{\phi_*^{-1}(x)}}\in\CF(\F)$. If $\phi$
is also representable, define $\LCF^\stk(\phi)f\!=\!\LCF^\na
(\phi)(m_\phi\!\cdot\!f)$ as in~\eq{ca5eq2}.
\label{ca5def3}
\end{dfn}

Suppose $C\subseteq\G(\K)$ is constructible, and let
$B=\phi_*^{-1}(C)$. Then $B$ is constructible as $\phi$
is of finite type, by Definition \ref{ca5def2}. Write
$\de_B,\de_C$ for the characteristic functions of $B,C$.
Then $f\cdot\de_B\in\CF(\F)$, and it follows easily that
\begin{equation*}
\bigl(\LCF^\na(\phi)f\bigr)\de_C\!=\!\CF^\na(\phi_*)(f\cdot\de_B)
\;\>\text{and}\;\>
\bigl(\LCF^\stk(\phi)f\bigr)\de_C\!=\!\CF^\stk(\phi)(f\cdot\de_B).
\end{equation*}
Therefore $\bigl(\LCF^\na(\phi)f\bigr)\vert_C$, $\bigl(\LCF^\stk
(\phi)f\bigr)\vert_C$ are constructible by Theorems \ref{ca4thm}
and \ref{ca5thm2}, for any constructible $C\subseteq\G(\K)$.
Hence $\LCF^\na(\phi)f$, $\LCF^\stk(\phi)f$ are {\it locally
constructible}, and $\LCF^\na(\phi),\LCF^\stk(\phi)$ are
linear maps $\LCF(\F)\ra\LCF(\G)$. From Theorems \ref{ca4thm}
and \ref{ca5thm2} we deduce:

\begin{thm} Let\/ $\F,\G,\H$ be algebraic $\K$-stacks with affine
geometric stabilizers, and\/ $\phi:\F\!\ra\!\G$, $\psi:\G\!\ra\!\H$
finite type $1$-morphisms. Then so is $\psi\ci\phi$, and\/
$\LCF^\na(\psi\ci\phi)\!=\!\LCF^\na(\psi)\ci\LCF^\na(\phi)$.
If\/ $\phi,\psi$ are representable then so is $\psi\ci\phi$,
and\/~$\LCF^\stk(\psi\ci\phi)\!=\!\LCF^\stk(\psi)\ci\LCF^\stk(\phi)$.
\label{ca5thm4}
\end{thm}

The locally constructible analogue of Theorem \ref{ca5thm3} also holds.


\begin{thebibliography}{99}

\bibitem{AtSe} M. Atiyah and G. Segal, {\it On equivariant Euler
characteristics}, J. Geom. Phys. 6 (1989), 671--677.

\bibitem{DHVW} L. Dixon, J.A. Harvey, C. Vafa and E. Witten, {\it
Strings on orbifolds}, Nucl. Phys. B261 (1985), 678--686.

\bibitem{Gome} T.L. G\'omez, {\it Algebraic stacks}, Proc. Indian Acad.
Sci. Math. Sci. 111 (2001), 1--31. math.AG/9911199.

\bibitem{Grot} A. Grothendieck and others, {\it S\'eminaire de
G\'eom\'etrie Alg\'ebrique}. Part 4: {\it Th\'eorie des topos et
cohomologie \'etale des sch\'emas}, Lecture Notes in Math. 269,
270 and 305, Springer, Heidelberg, 1972-3. Part $4\frac{1}{2}$:
{\it Cohomologie \'etale} (by P. Deligne), Lecture Notes in
Math. 569, Springer, Heidelberg, 1977. Part 5: {\it Cohomologie
$l$-adique et fonctions $L$}, Lecture Notes in Math. 589,
Springer, Heidelberg, 1977.

\bibitem{Hart} R. Hartshorne, {\it Algebraic Geometry}, Graduate
Texts in Math. 52, Springer-Verlag, New York, 1977.

\bibitem{Illu} L. Illusie, {\it Th\'eorie de Brauer et
caract\'eristique d'Euler--Poincar\'e d'apr\`es P. Deligne},
Exp. VIII in Ast\'erisque 82-83, 1981.

\bibitem{Joyc1} D.D. Joyce, {\it Motivic invariants of Artin
stacks and `stack functions'}, math.AG/0509722, version 2, 2006.

\bibitem{Joyc2} D.D. Joyce, {\it Configurations in abelian
categories. I. Basic properties and moduli stacks}, math.AG/0312190,
version 5, 2005.

\bibitem{Joyc3} D.D. Joyce, {\it Configurations in abelian
categories. II. Ringel--Hall algebras}, math.AG/0503029, version 3,
2006.

\bibitem{Joyc4} D.D. Joyce, {\it Configurations in abelian
categories. III. Stability conditions and identities},
math.AG/0410267, version 4, 2006.

\bibitem{Joyc5} D.D. Joyce, {\it Configurations in abelian
categories. IV. Invariants and changing stability conditions},
math.AG/0410268, version 4, 2006.

\bibitem{KaLa} N.M. Katz and G. Laumon, {\it Transformation de
Fourier et majoration des sommes exponentielles}, Publ. Math. I.H.E.S.
62 (1985), 145--202.

\bibitem{Kenn} G. Kennedy, {\it MacPherson's Chern classes of
singular algebraic varieties}, Communications in Algebra 18 (1990),
2821--2839.

\bibitem{Kres} A. Kresch, {\it Cycle groups for Artin stacks},
Invent. math. 138 (1999), 495--536. math.AG/9810166.

\bibitem{LaMo} G. Laumon and L. Moret-Bailly, {\it Champs alg\'ebriques},
Ergeb. der Math. und ihrer Grenzgebiete 39, Springer-Verlag, Berlin, 2000.

\bibitem{MacP} R.D. MacPherson, {\it Chern classes for singular
algebraic varieties}, Ann. Math. 100 (1974), 423--432.

\bibitem{Mark} D. Marker, {\it Model Theory: An Introduction},
Graduate Texts in Math. 217, Springer, New York, 2002.

\bibitem{Miln} J.S. Milne, {\it \'Etale Cohomology}, Princeton
University Press, Princeton, 1980.

\bibitem{Mumf} D. Mumford, {\it The Red Book of Varieties and
Schemes}, second edition, Lecture Notes in Math. 1358,
Springer-Verlag, Berlin, 1999.

\bibitem{Roan} S.-S. Roan, {\it Orbifold Euler Characteristic}.
In B. Greene and S.-T. Yau, editors, {\it Mirror Symmetry II},
pages 129--140. A.M.S. and International Press, 1997.

\bibitem{Rose} M. Rosenlicht, {\it A Remark on Quotient Spaces},
Annaes da Academia Brasileira de Ci\^encas 35 (1963), 487--489.

\bibitem{Viro} O.Y. Viro, {\it Some integral calculus based on
Euler characteristic}. In O.Y. Viro, editor, {\it Topology and
Geometry -- Rohlin seminar}, Lecture Notes in Math. 1346, pages
127--138. Springer-Verlag, Berlin, 1988.

\end{thebibliography}
\end{document}